\documentclass[11pt]{amsart}
\usepackage{latexsym,amsmath,amssymb,amsthm}
\theoremstyle{plain}
\newtheorem{thm}{Theorem}[section]

\newtheorem{cor}[thm]{Corollary}
\newtheorem{prop}[thm]{Proposition}
\newtheorem{rem}[thm]{Remark}

\DeclareMathOperator{\tr}{tr} \DeclareMathOperator{\Tr}{Tr}

\theoremstyle{definition}

\renewcommand{\Re}{\mathrm{Re}}
\newcommand{\Det}{\mathrm{Det}}
\newcommand{\G}{\Gamma}

\newcommand{\Spin}{\mathrm{Spin}}

\newcommand{\ski}{\sum_{i=1}^\kappa}

\newcommand{\Dd}{\mathcal D}

\newcommand{\Id}{\mathrm{Id}}

\newcommand{\bs}{\backslash}

\newcommand{\dsi}{d(\sigma_+)}

\newcommand{\bC}{\mathbb{C}}

\begin{document}
\baselineskip.550cm
\title[Eta invariant, regularized determinant ]
{Eta invariants and regularized determinants for odd dimensional
hyperbolic manifolds with cusps}

\author{Jinsung Park}

\address{Mathematisches Institut\\ Universit\"at Bonn\\
Beringstra{\ss}e 1\\ D-53115 Bonn\\ Germany}
\email{jpark@math.uni-bonn.de}


\maketitle

\begin{abstract}

We study eta invariants of Dirac operators and regularized
determinants of Dirac Laplacians over hyperbolic manifolds with
cusps and their relations with Selberg zeta functions.  Using the
Selberg trace formula and a detailed analysis of the unipotent
orbital integral, we show that the eta and zeta functions defined
by the relative traces are regular at the origin so that we can
define the eta invariant and the regularized determinant. We also
show that the Selberg zeta function of odd type has a meromorphic
extension over $\mathbb{C}$, prove a relation of the eta invariant
and a certain value of the Selberg zeta function of odd type, and
derive a corresponding functional equation. These results
generalize the earlier work of John Millson (see \cite{M}) to
hyperbolic manifolds with cusps. We also prove that the Selberg
zeta function of even type has a meromorphic extension over
$\mathbb{C}$, relate it to the regularized determinant, and obtain
a corresponding functional equation.
\end{abstract}

\section{Introduction}\label{s0}

In the seminal paper \cite{M}, Millson derived a relation of the
eta invariant of the odd signature operator and a certain value of
Selberg zeta function of odd type for compact hyperbolic manifolds
of dimension $(4n - 1)$. To prove this, Millson used the Selberg
trace formula, which relates the spectral data to the geometric
data, applied to a test function defined by the odd heat kernel of
the odd signature operator. In \cite{F}, the corresponding work
was done for the analytic torsion and the Ruelle zeta function for
odd dimensional hyperbolic manifolds. These results have been
generalized to the case of compact locally symmetric spaces of
higher rank in \cite{MS1}, \cite{MS2} .

It would be an interesting problem to extend the aforementioned
results to noncompact locally symmetric spaces with finite
volumes. However, we encounter several serious difficulties when
we discuss the extension of those results to noncompact locally
symmetric spaces. First of all, the heat operator of the Laplacian
is not of trace class. Hence we cannot use the trace of the heat
operator as in the compact case. In \cite{Mu3}, \cite{Mu4},
M\"uller introduced the relative trace to overcome this kind of
difficulty and defined corresponding relative eta invariants and
relative regularized determinants. In this paper, we follow
M\"uller's approach to study the (relative) spectral invariants
for noncompact hyperbolic manifolds with finite volumes. If the
continuous spectrum of the operator has a gap near zero, the
relative trace behaves as the usual trace in the compact case.
This is the case of the Laplacians acting on functions. The
regularized determinants of these operators and their relations
with Selberg zeta functions have been extensively studied for two
and three dimensional manifolds with cusps  in \cite{E1},
\cite{E2}, \cite{K1}, \cite{K2}, \cite{Mu2}. However, if the
continuous spectrum of the operator has no gap near zero, we
need to know the large time behavior of the relative trace, so its
relation with the scattering theory. In this paper, we study the
spectral invariants of Dirac operators and Dirac Laplacians for
odd dimensional hyperbolic manifolds with cusps, whose continuous
spectrums reach zero.

Our approach is to use the Selberg trace formula applied to test
functions defined by the heat kernels. To do so, we analyze the
corresponding geometric side of the Selberg trace formula, in
particular, the unipotent orbital integral. A detailed analysis of
these terms enables us to show that the Selberg zeta functions of
odd/even type have meromorphic extensions over $\mathbb{C}$. These
results can be considered as generalizations of an old result of
Gangolli and Warner in \cite{GW} to the case of nontrivial locally
homogeneous vector bundles over noncompact locally symmetric
spaces.

We explain our result more precisely. Let $X=\G\backslash
\Spin(2n+1,1)/\Spin(2n+1)$ be a $(2n+1)$-dimensional hyperbolic
manifold with cusps. Here $\G$ is a discrete subgroup of
$G=\Spin(2n+1,1)$ with finite co-volume. Throughout this paper, we
also assume that the group generated by the eigenvalues of $\G$
contains no root of unity. We now consider the Dirac operator $\Dd$
acting on $L^2(X,E)$. Here the spinor bundle $E$ over $X$ is a
locally homogeneous vector bundle defined by the spin representation
$\tau_{n}$ of the maximal compact subgroup $\Spin(2n+1)$ of
$\Spin(2n+1,1)$. Let us observe that the restriction of $\tau_n$ to
$\Spin(2n)\subset \Spin(2n+1)$ has the decomposition
$\sigma_+\oplus\sigma_-$ where $\sigma_{\pm}$ denotes the half spin
representation of $\Spin(2n)$.  Let $\mathcal{K}_t$ be the family of
functions over $G=\Spin(2n+1,1)$ given by taking the local trace of
the integral kernel of $e^{-t\widetilde{D}^2}$ or $\widetilde{D}
e^{-t\widetilde{D}^2}$ where $\widetilde{D}$ is the lifting of $\Dd$
over the universal covering space of $X$. Now the Selberg trace
formula applied to $\mathcal{K}_t$ has the following form,
\begin{align*}
&\sum_{\sigma=\sigma_{\pm}}\sum_{\lambda_k\in \sigma_p^{\pm}}
\widehat{\mathcal{K}}_{t}(\sigma,i\lambda_k)
-\frac{1}{4\pi}\int^{\infty}_{-\infty} \Tr(
C_{\Gamma}(\sigma_{+},-i\lambda)
C_{\Gamma}'(\sigma_{+},i\lambda)\pi_{\Gamma}(\sigma_{+},i\lambda)
({\mathcal{K}}_{t}) )
\ d\lambda \\
=&\quad
I_{\G}(\mathcal{K}_t)+H_{\G}(\mathcal{K}_t)+U_{\G}(\mathcal{K}_t)
\end{align*}
where $\sigma_p:=\sigma^+_p\cup \sigma^{-}_p$ gives the point
spectrum of $\Dd$ , $C_{\Gamma}(\sigma_{+},i\lambda)$ is the
intertwining operator (note that $C_{\Gamma}(\sigma_+,i\lambda)=
C_{\Gamma}(\sigma_-,i\lambda)$ since $\sigma_{+},\sigma_{-}$ are
unramified) and $I_{\G}(\cdot)$, $H_{\G}(\cdot)$, $U_{\G}(\cdot)$
are the identity, hyperbolic and unipotent orbital integrals,
respectively. Following \cite{Mu1}, \cite{Mu3}, \cite{Mu4}, we
define certain operators $\Dd_0(i)$ determined by $\Dd$ such that
$e^{-t\Dd^2}- \sum^\kappa_{i=1}e^{-t\Dd_0(i)^2}$ , $\Dd
e^{-t\Dd^2}- \sum_{i=1}^\kappa\Dd_0(i) e^{-t\Dd_0^2(i)}$ are trace
class operators on $L^2(X,E)$ where $\kappa$ is the number of the
cusps of $X$. We also show that the spectral side of the above
Selberg trace formula is equal to the relative trace $\Tr(
e^{-t\Dd^2}-\sum_{i=1}^\kappa e^{-t\Dd^2_0(i)})$ or $\Tr(\Dd
e^{-t\Dd^2}-\sum_{i=1}^\kappa \Dd_0(i) e^{-t\Dd^2_0(i)})$; for
example, if $\mathcal{K}_t$ is give by $e^{-t\widetilde{D}^2}$,
\begin{align*}
&\Tr(e^{-t\Dd^2}-\sum_{i=1}^\kappa e^{-t\Dd^2_0(i)})\\
 =& \sum_{\sigma=\sigma_{\pm}}\sum_{\lambda_k\in \sigma_p^{\pm}}\widehat{\mathcal{K}}_{t}(\sigma,i\lambda_k)
-\frac{1}{4\pi}\int^{\infty}_{-\infty}
\Tr(C_{\Gamma}(\sigma_{+},-i\lambda)
C_{\Gamma}'(\sigma_{+},i\lambda)\pi_{\Gamma}(\sigma_{+},i\lambda)({\mathcal{K}}_{t}))
\ d\lambda \\
=&\quad
I_{\G}(\mathcal{K}_t)+H_{\G}(\mathcal{K}_t)+U_{\G}(\mathcal{K}_t).
\end{align*}
A similar formula holds for $\Tr(\Dd
e^{-t\Dd^2}-\sum_{i=1}^\kappa\Dd_0(i)e^{-t\Dd^2_0(i)})$ with the
corresponding kernel function $\mathcal{K}_t$. The exact forms of
$I_{\G}(\mathcal{K}_t)$, $H_{\G}(\mathcal{K}_t)$ are well-known,
in particular, the hyperbolic orbital integrals
$H_{\G}(\mathcal{K}_t)$ provide us with the Selberg zeta function
of odd/even type. Hence, a main task in this paper is the analysis
of $U_{\G}(\mathcal{K}_t)$. To do so, we use a result of Hoffmann
in \cite{H} and perform explicit computations for the weighted
unipotent orbital integrals for our concerned cases. These
explicit computations constitute some of the main ingredients of
this paper. By these explicit computations, we can show that
$U_{\G}(\mathcal{K}_t)=0$ if $\mathcal{K}_t$ is determined by
$\widetilde{D} e^{-t\widetilde{D}^2}$.

Following M\"uller's approach in \cite{Mu2}, \cite{Mu3},
\cite{Mu4}, we define the eta function $\eta_{\Dd}(s)$ and the
zeta function $\zeta_{\Dd^2}(s)$ using the relative traces
$\Tr(e^{-t\Dd^2}-\sum^\kappa_{i=1}e^{-t\Dd^2_0(i)})$, $\Tr(\Dd
e^{-t\Dd^2}-\sum_{i=1}^\kappa\Dd_0(i) e^{-t\Dd^2_0(i)})$,
respectively (see \eqref{eta-zeta}, \eqref{def-eta},
\eqref{def-zeta} for the precise definitions of $\eta_{\Dd}(s)$,
$\zeta_{\Dd^2}(s)$ ). The continuous spectrum of $\Dd$ is the
whole real line, hence the large time contributions of the
relative traces can also give rise to poles of $\eta_{\Dd}(s)$,
$\zeta_{\Dd^2}(s)$. Therefore, we need to consider separately the
small time contributions and the large time contributions for the
meromorphic extensions of $\eta_{\Dd}(s)$, $\zeta_{\Dd^2}(s)$. We
use the analytic expansion of the intertwining operator along the
imaginary axis to get the large time contribution and we analyze
all the terms in the geometric side of the Selberg trace formula
for the small time contribution. We prove the following
meromorphic structures of the eta and zeta functions.

\begin{thm}\label{t:pole} The eta function $\eta_{\Dd}(z)$ and
the zeta function $\zeta_{\Dd^2}(z)$ have the meromorphic
structures,
\begin{align*}
&\Gamma((z+1)/2)\,\eta_{\Dd}(z)=\sum^{\infty}_{k=0}\frac{-2\gamma_k}{z-2k-2}+K(z),\\
&\Gamma(z)\,
\zeta_{\Dd^2}(z)=\sum^{\infty}_{k=-n}\frac{\beta_k}{z+k-\frac 12}
+\frac{\beta'_0}{(z-\frac 12)^2}-\frac{h}{z}
+\sum^{\infty}_{k=0}\frac{-\gamma'_k}{z-k-\frac 12}+H(z)
\end{align*}
where $\beta_k$, $\beta'_0$ are locally computable constants,
$\gamma_k$, $\gamma'_k$ are constants which are determined by the
intertwining operator $C_{\Gamma}(\sigma_{+}, i\lambda)$, $h$ is
the multiplicity of the zero eigenvalues of $\Dd$, and $K(z)$,
$H(z)$ are holomorphic functions . In particular, $\eta_{\Dd}(z)$
and $\zeta_{\Dd^2}(z)$ are regular at $z=0$.
\end{thm}

It follows from Theorem \ref{t:pole} that we can define the eta
invariant by
$$\eta(\Dd):=\eta_{\Dd}(0).
$$
Following \cite{M}, we use the Selberg trace formula to derive the
relation of the eta invariant $\eta(\Dd)$ and the value of the
Selberg zeta function of odd type $Z^o_H(s)$ at $s=n$ (see
\eqref{def-odd} for the precise definition of $Z^o_H(s)$).   Let
us remark that $Z^o_H(s)$ is defined a priori only for $\Re(s) \gg
0$ and the meromorphic extension of $Z^o_H(s)$ over $\bC$ is one
of the main results of this paper. As we mentioned above, we show
that all unipotent terms are vanishing in the Selberg trace
formula applied to the odd heat kernel function. In conclusion,
the equality of the relative trace and the orbital integrals for
the odd kernel function is exactly the same as in the case of
compact hyperbolic manifolds. Therefore we can expect the same
formula of $\eta(\Dd)$ and  $Z^o_H(s)$  as in compact hyperbolic
manifolds. However, in the corresponding functional equation, the
term determined by the intertwining operator appears. The
following theorem states our results for $\eta(\Dd)$ and
$Z^o_H(s)$.

\begin{thm}\label{t:eta0}
The Selberg zeta function of odd type $Z^o_H(s)$ has a meromorphic
extension over $\mathbb{C}$ with $s=n$ as a regular point and the
following equalities hold:
\begin{align*}
\eta(\Dd)&=\frac{1}{\pi i}\log Z^o_H(n),\\
Z^o_H(n+s)Z^o_H(n-s)&=\exp(2\pi i\eta(\Dd))\ \biggl( \frac{\det
C_+(s) C_-(0)}{\det C_-(s) C_+(0)}\biggr)^{-2^{n-1}} \qquad
\text{for} \quad s\in \mathbb{C} .\ \
\end{align*}
Here $C_{\pm}(s)$ are linear operators given by
$C_{\Gamma}(\sigma_+, s)=\begin{pmatrix} 0 & C_-(s)\\
                                         C_+(s) & 0 \end{pmatrix} $.

\end{thm}

We can define the zeta function $\zeta_{\Dd^2}(z,s)$ of the
shifted Dirac Laplacian $\Dd^2+s^2$ where $s$ is a positive real
number. Now the continuous spectrum of $\Dd^2+s^2$ does not reach
$0$ and the large time contribution does not create any poles of
$\zeta_{\Dd^2}(z,s)$. Therefore we can see that
$\zeta_{\Dd^2}(z,s)$ is regular at $z=0$ by Theorem \ref{t:pole}.
It follows that the regularized determinant
$$
\Det(\Dd^2,s):= \exp(-\zeta'_{\Dd^2}(0,s))
$$
is well-defined. We show that $\Det(\Dd^2,s)$ can be extended to a
meromorphic function of $s$ on $\mathbb{C}$.  Let us observe that
$\Det(\Dd^2,s)\neq \Det(\Dd^2,-s)$ as a meromorphic function over
$\bC$ (see Remark \ref{rem:even}). We use the Selberg trace
formula to prove a relation between $\Det(\Dd^2,s)$ and the
geometric data, which consists of the Selberg zeta function of
even type $Z^e_H(s)$ (see \eqref{def-even} for the precise
definition of $Z^e_H(s)$) and certain meromorphic functions over
$\bC$ derived from the identity and unipotent orbital integrals.
Let us also remark that $Z^e_H(s)$ is defined a priori only for
$\Re(s)\gg 0$ and its meromorphic extension over $\bC$ is one of
the main results of this paper. By a similar manner as in the
previous case, we also derive a functional equation for
$\Det(\Dd^2,s)$ and $Z^e_H(s)$ where the unipotent factor plays a
non-trivial role. The following theorem states our results for
$\Det(\Dd^2,s)$ and $Z^e_H(s)$.

\begin{thm}\label{t:det0}
The Selberg zeta function of even type $Z^e_H(s)$ has a
meromorphic extension over $\mathbb{C}$ and the following
equalities hold for any $s \in \mathbb{C}$:
\begin{align*}
\Det(\Dd^2,s)\ =&\ C\, Z^e_H(s+n) \Gamma(s+\frac
12)^{-2^{n}\kappa} \exp \big(4\pi\int^s_0
p(\sigma_{+},i\lambda)+P_U(i\lambda)\ d\lambda\big),\\
\Det(\Dd^2,s)^2 \ =& \ C^2\,  (\det
C_{\Gamma}(\sigma_+,s))^{-2^{n-1}}\ Z^e_H(n+s)Z^e_H(n-s)
\Bigr(\Gamma(s+\frac 12)\Gamma(-s+\frac 12)\Bigr)^{-2^{n}\kappa }
\notag
\end{align*}
where $C$ is a constant, $p(\sigma_{+},s)$ is Plancherel measure
for $ \sigma_+$, $P_U(s)$ is an even polynomial of degree
$(2n-4)$, and $\kappa$ denotes the number of the cusps of $X$ .
\end{thm}

If we compare Theorem \ref{t:eta0} with Theorem \ref{t:det0}, we
can see that there are no defect terms determined by the cusps in
the relation of $\eta(\Dd)$ and $Z^o_H(n)$. But the defect terms
appear in the relation of $\Det(\Dd^2,s)$ and $Z^e_H(s)$.

This paper has the following structure. In Section \ref{s1}, we
define a Dirac operator $\Dd$ for the locally homogeneous vector
bundle $E$, which is defined by the spin representation $\tau_n$
of maximal compact subgroup $\Spin(2n+1)$. We introduce an
operator $\Dd_0(i)$ naturally determined by $\Dd$, and using this
we define the relative trace. In Section \ref{s2}, we introduce
the Selberg trace formula for nontrivial locally homogeneous
vector bundles over noncompact locally symmetric spaces of rank
$1$. In Section \ref{s3}, we use result from \cite{H} to analyze
the unipotent terms. In Section \ref{s4}, we derive the relation
between the relative trace and the spectral side of the Selberg
trace formula. In Section \ref{s5}, we define $\eta_{\Dd}(s)$ and
$\zeta_{\Dd^2}(s)$ and prove Theorem \ref{t:pole}. In Sections
\ref{s6} and \ref{s7}, we prove Theorem \ref{t:eta0} and Theorem
\ref{t:det0} using results of Sections \ref{s2}, \ref{s3} and
\ref{s4}.

 {\bf Acknowledgments.} The author wants to express his gratitude
to Werner M\"uller and Werner Hoffmann for their helpful comments
on this paper. He also thanks Paul Loya, Morten Skarsholm Risager,
Masato Wakayama and Krzysztof Wojciechowski for several help
during the writing of this paper. Finally he thanks the anonymous
referee for pointing out many mistakes and giving several
comments, which improve this paper considerably. A part of this
work was done during the author's stay at ICTP and MPI. He wishes
to express his thanks to ICTP and MPI for their financial support
and hospitality.

\section{Dirac operators on odd dimensional hyperbolic manifolds
with cusps}\label{s1}

Let $G$ be a noncompact connected simple Lie group with finite
center and let $\frak{g}$ denote its Lie algebra. We assume that the
real rank of $G$ is one.  Let $\G$ be a discrete subgroup of $G$
such that the group generated by the eigenvalues of $\G$ contains no
root of unity. Let $K$ be a maximal compact subgroup of $G$,
$\theta$ the associated Cartan involution, $\frak{k}\oplus\frak{p}$
the associated Cartan decomposition, and $C(\ ,\ )$ the Killing form
of $\frak{g}$. Let $\frak{a}$ be a maximal abelian subalgebra of
$\frak{p}$; $\frak{a}$ is one dimensional since the rank of $G$ is
$1$. Let $\Phi$ be the set of roots of $(\frak{g}, \frak{a})$.
Choose an order for $\Phi$ and let $\frak{g}=\frak{g}_0\oplus
\sum_{\lambda\in\Phi^+}\frak{g}_{\lambda}\oplus
\sum_{\lambda\in\Phi^+}\frak{g}_{-\lambda}$ be the root space
decomposition where $\frak{g}_0$ is the centralizer of $\frak{a}$.
Let $\alpha_1$ be the unique simple positive root; then
$\Phi^+=\{\alpha_1, 2\alpha_1\}$. Let
$\frak{n}_{\alpha_1}:=\frak{g}_{\alpha_1}$ and
$\frak{n}_{2\alpha_1}:=\frak{g}_{2\alpha_1}$. This last space may be
$0$. If $N=\exp(\frak{n}_{\alpha_1}\oplus \frak{n}_{2\alpha_1})$ and
$A=\exp(\frak{a})$, then $G=NAK$ is an Iwasawa decomposition. Let
$P_0=NAM$ (with $M$ the centralizer of $A$ in $K$) be the associated
minimal parabolic. The $G$-conjugates of $P_0$ are the proper
parabolic subgroups of $G$. A parabolic subgroup $P$ is called
$\G$-cuspidal if $\G\cap N(P)\bs N(P)$ is compact. Here $N(P)$ is
the unipotent radical of $P$ which, we may assume, is $G$-conjugate
to $N$. Let $P_{\G}=\{P_1,\cdots, P_{\kappa}\}$ be a complete set of
$\G$-conjugacy classes of $\G$-cuspidal subgroups of $G$.

From now on, we assume that $G=\Spin(2n+1,1)$ and $K=\Spin(2n+1)$.
Let $\widetilde X$ denote the noncompact symmetric space given by $
\Spin(2n+1,1)/\Spin(2n+1)$.  The Killing form $C( \ , \ )$ provides
us with an invariant metric on $G/K$ by $\langle Y ,Z
\rangle=\frac{1}{4n}C( Y , Z )$ for $Y,Z\in \frak{p}$, which gives
us constant curvature $(-1)$. We use the spin representation $(
\tau_{n},V_{\tau_n} ) $ of $\Spin(2n+1)$ to define a homogeneous
vector bundle $\widetilde {E}$ over
$\widetilde{X}=\Spin(2n+1,1)/\Spin(2n+1)$ by $\widetilde{E}=G\times
V_{\tau_{n}}/\sim $  where
$$(g,v)\sim (g', v') \quad \ \ \text{if} \ \ \ (g', v')=(gk, \tau_{n}(k^{-1})v)$$
for $g, g' \in G,\  k\in K$. We denote such equivalence classes by
$[g,v]$. Note that $\widetilde{E}$ admits a left $G$ action
defined by $g_0[g,v]=[g_0g,v]$. If we restrict the spin
representation $\tau_{n}$ of $K=\Spin(2n+1)$ to $M=\Spin(2n)$,
then $\tau_n$ decomposes into two half spin representations
$\sigma_{+}, \sigma_{-}$ of $\Spin(2n)$ . Two representations
$(\sigma_{+}, H_{\sigma_+}), (\sigma_{-}, H_{\sigma_-})$ are
unramified and $w\sigma_{+}=\sigma_{-},w\sigma_{-}=\sigma_{+}$ for
the nontrivial element $w\in W(A)=M^*/M$ where $M^*$ is the
normalizer of $A$ in $K$. We define the Dirac operator
$\widetilde{D}:C^{\infty}(\widetilde X, \widetilde {E})\ \to \
C^{\infty}(\widetilde X, \widetilde {E})$ by
$$\widetilde{D}=\sum^{2n+1}_{i=1} c(X_i)  \nabla_{X_i}$$
where $\{X_i:1\leq i \leq 2n+1\}$ is a left invariant orthonormal
frame such that $H:=X_{2n+1}$ at $eK$ spans $\frak{a}$, and
$\nabla$ is the Levi-Civita connection on $\widetilde{E}$.

We consider a locally symmetric space given by $X=\G\bs
\Spin(2n+1,1)/\Spin(2n+1)$ where $\G$ is a discrete subgroup of
$\Spin(2n+1,1)$ with unipotent elements satisfying the condition
in the introduction. As a consequence of this assumption, $\G$ is
torsion free and $\G\cap P=\G\cap N(P)$ so that $\G \cap P\bs
N(P)= \G\cap N(P) \bs N(P) $. Then $X$ is a $(2n+1)$-dimensional
hyperbolic manifold with cusps. We denote by $E$ the quotient
space $\G \bs \widetilde{E}$. This is a locally homogeneous vector
bundle over $X$. Moreover the Dirac operator $\widetilde{D}$ can
be pushed down to $X$. We denote this operator by $D$ and its
unique self adjoint extension on $L^2(X,E)$ by $\Dd$. The
hyperbolic manifold $X$ endowed with the metric $\langle\ , \
\rangle$ has a following decomposition,
\begin{equation}\label{decompX}
X=X_0\cup W_1 \cup \cdots \cup W_\kappa
\end{equation}
where $X_0$ is a compact manifold with boundaries and $W_i$,
$i=1,\cdots, \kappa$ are ends of $X$.
(In general, $X_0$ and $W_i$'s may have the nonempty intersections each other.)
For each end $W_i$ which we call as \emph{cusp}, $W_i \cong
[0,\infty)\times N_i$ and $N_i$ can be identified with the flat
torus $T^{2n}$ with the metric $dn^2$, and the restriction of
$\langle \ , \ \rangle$ to $W_i$ has the form $ dg^2(r,x)=dr^2 +
e^{-2r}dn^2(x) $ for $(r,x)\in [0,\infty)\times N_i$.

Now we have the following expression for $D$ over the cusp $W_i$ ,
\begin{equation}\label{e:dirac}
D=c(H)(\nabla_H + B - n\, \Id)
\end{equation}
where $B=\sum^{2n}_{i=1}c(X_i)c(H)\nabla^N_{X_i}$ with Levi-Civita
connection $\nabla^N$ over $E|_N$ (see (7) in \cite{B}). Note that
$B$ has the property $ c(H) B=-B c(H)$ and $E_0:=\ker(B)$ can be
identified with $V_{\tau_n}$. Let $0<\mu_1 \leq \mu_2 \leq \mu_3
\leq \cdots \to\infty$ be positive eigenvalues of $B$, each
eigenvalue repeated according to its multiplicity with
corresponding eigensections $\psi_j$. We denote by $E_{\mu_j}$ the
eigenspace corresponding to $\mu_j$. We decompose
$$
L^2(\mathbb{R}^+\times N_i, E\vert_{\mathbb{R}^{+}\times N_i},
dg^2(r,x))
$$
into
$$
L^2(\mathbb{R}^+, E_0, e^{-2nr}dr)\oplus \oplus_{j=1}^{\infty}
L^2(\mathbb{R}^+,E_{\mu_j}\oplus c(H)E_{\mu_j}, e^{-2nr}dr)
$$
where $\mathbb{R}^+$ denotes $[0,\infty)$. Note that $E_0\cong
V_{\tau_n}$ is a symplectic vector space with a symplectic
structure $\langle c(H)\ \ , \ \rangle$. We fix a Lagrangian
subspace $L$ of $E_0\cong V_{\tau_n}$ such that $ E_0= L\oplus
c(H)L $. Then the map
$$
\alpha_j\psi_j + \alpha_{-j}c(H)\psi_j \rightarrow \  e^{-nr}
\left(\begin{matrix}
\alpha_j\\
\alpha_{-j}
\end{matrix}\right )
$$
gives a unitary equivalence
$$
L^2(\mathbb{R}^+, E_0, e^{-2nr}dr)\oplus \oplus_{j=1}^{\infty}
L^2(\mathbb{R}^+,E_{\mu_j}\oplus c(H)E_{\mu_j}, e^{-2nr}dr)
$$
$$
\cong  \quad \oplus^{d(\sigma_{+})}_{j=1}
L^2(\mathbb{R}^+,\mathbb{C}^2,dr)\oplus\oplus^{\infty}_{j=1}L^2(\mathbb{R}^+,\mathbb{C}^2,dr)
$$
under which the Dirac operator $D$ over the cusp $W_i$ decomposes
as follows
\begin{equation}\label{e:d0}
D \rightarrow  \ \ \oplus^{d(\sigma_{+})}_{i=1} D_0 \oplus \oplus
^{\infty}_{j=1} D_{\mathbb{\mu}_j}
\end{equation}
where
$$
D_{\mu}= \left( \begin{matrix}
0 & -\frac{d}{dr} + \mu \\
 \frac{d}{dr} + \mu  & 0
\end{matrix} \right).
$$
Hence the Dirac Laplacian $D^2$ is transformed into
$\oplus_{\mu\in \mathrm{spec}(B)} D^2_{\mu} $ where
$D^2_{\mu}=-\frac{d^2}{dr^2} +\mu^2$ on
$L^2(\mathbb{R}^{+},\mathbb{C},dr)$.

We now consider the operator
$$
c(H)(\frac{d}{dr}-n\, \Id):C_0^{\infty}(\mathbb{R}^{+},E_0) \to
C_0^{\infty}(\mathbb{R}^{+},E_0 )
$$
whose $L^2$-extension ( with respect to $e^{-2nr}dr$ ) is
transformed to $d(\sigma_{+})$-copies of $D_0$ in (\ref{e:d0}).
Now, we put
\[
C_0^{\infty}(\mathbb{R}^{+},E_0,L) :=\{\ \phi\in
C_0^{\infty}(\mathbb{R}^{+},E_0)\ | \ \phi(0)\in L\ \} , \] then
the following operator
\begin{equation}\label{e:defD0}
c(H)(\frac{d}{dr}-n\, \Id) : C_0^{\infty}(\mathbb{R}^{+},E_0,L) \
\to \ L^2(\mathbb{R}^{+},E_0, e^{-2nr}dr)
\end{equation}
is essentially self adjoint.  By the natural embedding of
$\mathbb{R}^+$ into the geodesic rays in $W_i\subset X$, we can
regard $L^2(\mathbb{R}^+, E_0, e^{-2nr}dr )$ as a subspace of
$L^2(X,E)$. The operator $c(H)(\frac{d}{dr}-n\, \Id)$ in
\eqref{e:defD0} can be extended to the self adjoint operator on
$L^2(X,E)$ by the zero map over the orthogonal complement of this
subspace in $L^2(X,E)$. For each $W_i$, in this way we obtain the
operator $\Dd_0(i)$, $i=1,\cdots,\kappa$. We can see that each
$\Dd_0(i)$ has no point spectrum.  Now we have

\begin{prop} The differences $\big(e^{-t\Dd^2}-\sum_{i=1}^\kappa e^{-t\Dd^2_0(i)}\big)$
, $\big(\Dd e^{-t\Dd^2}- \sum_{i=1}^\kappa
\Dd_0(i)e^{-t\Dd^2_0(i)}\big)$ are trace class operators on
$L^2(X,E)$ for $t>0$.
\end{prop}

\begin{proof} First choose a smooth function $g_s$ such that
$0<g_s\le1$, $g_s(x)=1$ for $x\in X_0$ and $g_s((r,\cdot))=0$ for
$(r,\cdot)\in W_i$ with $r\ge s $ for $1\leq i \leq \kappa$.
Denote by $U_{g_s}$, $U_{1-g_s}$ the operators in $L^2(X,E)$
defined by multiplication by $g_s,1-g_s$ respectively. Then, for
$s\gg 0$, $U_{g_s}(\sum_{i=1}^\kappa e^{-t\Dd^2_0(i)})$ is of
trace class and the support of $U_{1-g_s}(\sum_{i=1}^\kappa
e^{-t\Dd^2_0(i)})$ consists of disjoint $\kappa$-components lying
in $W_i$'s. So, for our proof we may assume that the supports of $
e^{-t\Dd^2_0(i)}$, $W_i$'s are disjoint each other. Now pick $f\in
C^{\infty}(X)$ such that $0 < f \leq 1$, $f(x)=1$ for $x\in X_0$
and $f((r,\cdot))=e^{-\frac{r}{4}}$ for $(r,\cdot)\in W_i$ with
$r\gg 0 $ for $1\leq i \leq \kappa$. Denote by $U_f$ the operator
in $L^2(X,E)$ defined by multiplication by $f$. Then we may write
\begin{align*}
e^{-t\Dd^2}-\sum_{i=1}^\kappa e^{-t\Dd^2_0(i)} \,  =& \, (
e^{-\frac{t}{2}\Dd^2} - \ski e^{-\frac{t}{2}\Dd^2_0(i)} )\circ
U_f^{-1}\circ U_f\circ e^{-\frac{t}{2}\Dd^2}\\
&  + \ski e^{-\frac{t}{2}\Dd^2_0(i)}\circ U_f\circ U_f^{-1}\circ (
e^{-\frac{t}{2}\Dd^2} - \ski e^{-\frac{t}{2}\Dd^2_0(i)}).
\end{align*}
The heat kernel estimates for $e^{-t\Dd^2}$ and $\ski
e^{-t\Dd^2_0(i)}$ show that $( e^{-\frac{t}{2}\Dd^2} - \ski
e^{-\frac{t}{2}\Dd^2_0(i)} )\circ U_f^{-1}$, $ U_f\circ
e^{-\frac{t}{2}\Dd^2}$, $\ski e^{-\frac{t}{2}\Dd^2_0(i)}\circ U_f$
and $U_f^{-1}\circ ( e^{-\frac{t}{2}\Dd^2} - \ski
e^{-\frac{t}{2}\Dd^2_0(i)})$ are Hilbert-Schmidt operators. The
composition of Hilbert-Schmidt operators is of trace class, hence
$e^{-t\Dd^2}-\ski e^{-t\Dd^2_0(i)}$ is a trace class operator. The
remaining case is proved in the same way.
\end{proof}

\section{The Selberg trace formula}\label{s2}

 Let $R_{\G}$ be the right regular representation of $G$ on
$L^2(\G\bs G)$ and $f$ be a right $K$-finite function in the
Harish-Chandra $L^p$-Schwartz space $\mathcal{C}^p(G)$ where $0 <
p < 1$. The trace of the restriction $R_{\G}(f)$  to the discrete
part $L^2_{d}(\G\bs G)$ of $L^2(\G\bs G)$ can be written as
$$
\Tr(R_{\G}(f)\vert_{L^2_{d}(\G\bs
G)})=I_{\G}(f)+H_{\G}(f)+U_{\G}(f)+S_{\G}(f)+T_{\G}(f) .
$$
Here
\begin{align*}
I_{\G}(f)&=\text{Vol}(\G\bs G)f(1),\\
H_{\G}(f)&=\sum_{\{\gamma:\text{hyperbolic}\}}\text{Vol}(\G_{\gamma}\bs
G_{\gamma}) \int_{G_{\gamma}\bs G} f(x^{-1}\gamma x) dx,\\
U_{\G}(f)&=\frac{1}{\vert\alpha_1\vert}\{C_1(\G)T_1(f)+C_2(\G)T_2(f)+C'_1(\G)T'_1(f)\},\\
S_{\G}(f)&=\frac{1}{4\pi i
}\sum_{\sigma\in\hat{M}}\int_{\Re(\nu)=0}
\Tr(C_{\G}(\sigma,\nu)^{-1}\frac{d}{d\nu}C_{\G}(\sigma,\nu)\pi_{\G}(\sigma,\nu)(f))d\nu,\\
T_{\G}(f)&=-\frac{1}{4}\sum_{\sigma\in\hat{M}}\Tr(C_{\G}(\sigma,0)\pi_{\G}(\sigma,0)(f))
\end{align*}
where $\G_{\gamma}, G_{\gamma}$ are centralizers of $\gamma$ in
$\G, G$ respectively, the constants $C_1(\G)$, $C_2(\G)$,
$C'_1(\G)$ and tempered distributions $T_1,T_2,T'_1$ will be
discussed in the next section, and
$$
C_{\G}(\sigma,\nu):\mathcal{H}_{\G}(\sigma,\nu) \to
\mathcal{H}_{\G}(\sigma,-\nu)$$ denotes the intertwining operator
(see p.166--167 in \cite{BM}, p.9--10 in \cite{SW} for the precise
definitions of $C_{\G}(\sigma, \nu),
\mathcal{H}_{\G}(\sigma,\nu)$). We refer to \cite{BM}, \cite{OW},
\cite{SW}, \cite{War} for detailed expositions of the Selberg
trace formula. The various invariant measures that we use are
normalized as in \cite{OW}. More precisely, let $d_K$ be the Haar
measure on $K$ which assigns to $K$ total volume one, let $d_A$,
$d_N$ be the exponentiation of the normalized Lebesgue measure on
the Lie algebra $\frak{a}$, $\frak{n}$ of $A$, $N$, respectively,
relative to the Euclidean structure associated with the Killing
form. Then the Haar measure $d_G$ is determined by
\[
d_G(x)=a^{2\rho}\, d_N(n)\, d_A(a)\, d_K(k) \quad \text{for} \quad
x=nak
\]
where $\rho$ is the sum of the positive roots of
$(\frak{g},\frak{a})$ divided by $2$.

The kernel $\widetilde{K}_t(x:y)$ of the integral operator
$e^{-t\widetilde{D}^2}$ (or $\widetilde{D} e^{-t\widetilde{D}^2}$)
over $G/K$ is a section of $\widetilde{E}\boxtimes
\widetilde{E}^{*}$, the external tensor product of $\widetilde{E}$
and $\widetilde{E}^{*}$ over $G\times G$. The bundle $E$ is a
trivial bundle, hence the kernel $\widetilde{K}_t(x:y)$ is an
element of $(C^{\infty}(G\times G)\otimes
\text{End}(V_{\tau_{n}}))^K$, which consists of endomorphism
valued functions on $G\times G$ invariant under $K$. It follows
that there exists a function $ \widetilde{K}_t(x):G\to
\text{End}(V_{\tau_{n}}) $ such that
$$\widetilde{K}_t(x:y)=\widetilde{K}_t(x^{-1}y) \quad
\text{and} \quad
\widetilde{K}_t(k_1^{-1}xk_2)=\tau_{n}(k_1^{-1})\widetilde{K}_t(x)\tau_{n}(k_2)
$$
for $x,y \in G, k_1,k_2\in K$.  We define the local trace of
$\widetilde{K}_t(x)$ to be the scalar function on $G$ given by
${K}_t(x):=\tr(\widetilde{K}_t(x))$. We denote by
$\widehat{K}_t(\sigma,\nu)$ the Fourier transform of ${K}_t$ for
the unitary principal representation $\pi_{\sigma,\nu}$ of $G$,
that is,
$$
\widehat{K}_t(\sigma,\nu):=\Tr(\pi_{\sigma,\nu}(K_t)).
$$
From now on, we will denote by $K^e_t$, $K^o_t$ the scalar
functions corresponding to $e^{-t\widetilde{D}^2}$, $\widetilde{D}
e^{-t\widetilde{D}^2}$. By (4.5) in \cite{MS1}, we have
\begin{prop}\label{p:frie} For $\lambda \in \mathbb{R}$, we have
$$
\widehat{K}^e_t(\sigma_{\pm},i\lambda)=e^{-t\lambda^2}, \quad
\widehat{K}^o_t(\sigma_+,i\lambda)=\lambda e^{-t\lambda^2}, \quad
\widehat{K}^o_t(\sigma_{-},i\lambda)=-\lambda e^{-t\lambda^2}.
$$
\end{prop}

Since $\sigma_{\pm}$ is unramified, the intertwining operator $
C_{\G}(\sigma_+, \nu) = C_{\G}(\sigma_{-}, \nu) $ acting on
$\mathcal{H}_{\G}(\sigma_{+},\nu) =
\mathcal{H}_{\G}(\sigma_{-},\nu) $ switches the subspaces induced
by the representations $\sigma_+, \sigma_-$. Hence,
$C_{\G}(\sigma_+, \nu)$ takes the form
$$
C_{\G}(\sigma_+, \nu) = \left (
\begin{matrix} 0 & C_-(\nu) \\
C_+(\nu) & 0
\end{matrix}\right )
$$
with respect to the decomposition of
$\mathcal{H}_{\G}(\sigma_+,\nu)$.  Therefore, we have
\begin{equation*}
C_{\G}(\sigma_+,\nu)^{-1} C'_{\G}(\sigma_+,\nu)
=\left(\begin{matrix}  C_-(-\nu)C'_+(\nu) & 0 \\
0 & C_+(-\nu)C'_-(\nu)
\end{matrix}\right)
\end{equation*}
and
$$
T_{\G}(\cdot)=-\frac{1}{4}\Tr(C_{\G}(\sigma_+,0)\pi_{\G}(\sigma_+,0)(\cdot))=0.
$$
Since $\tau_n|_M=\sigma_+\oplus \sigma_-$, the Selberg trace
formulas applied to $K^e_t$, $K^o_t$ are given by
\begin{align}
\label{e:etrace}&\sum_{\sigma=\sigma_{\pm}}\sum_{\lambda_k\in
\sigma_p^{\pm}} \widehat{{K}}^e_{t}(\sigma,i\lambda_k)
-\frac{1}{4\pi}\int^{\infty}_{-\infty}
\Tr(C_{\Gamma}(\sigma_+,-i\lambda)
C_{\Gamma}'(\sigma_+,i\lambda)\pi_{\Gamma}(\sigma_+,i\lambda)
({K}^e_{t}))
\ d\lambda \\
=&\quad \sum_{\lambda_k\in \sigma_p}
e^{-t\lambda_k^2}-\frac{\dsi}{4\pi} \int^{\infty}_{-\infty}
e^{-t\lambda^2}\tr (C_{\G}(\sigma_{+},-i\lambda)
C_{\G}'(\sigma_{+},i\lambda))
\ d\lambda \notag \\
=&\quad I_{\G}(K^e_t)+H_{\G}(K^e_t)+U_{\G}(K^e_t) ,\notag
\end{align}

\begin{align}\label{e:otrace}
&\sum_{\sigma=\sigma_{\pm}}\sum_{\lambda_k\in \sigma_p^{\pm}}
\widehat{{K}}^o_{t}(\sigma,i\lambda_k)
-\frac{1}{4\pi}\int^{\infty}_{-\infty}
\Tr(C_{\Gamma}(\sigma_+,-i\lambda)
C_{\Gamma}'(\sigma_+,i\lambda)\pi_{\Gamma}(\sigma_+,i\lambda)(K^o_t))
\ d\lambda \\
=&\quad \sum_{\lambda_k\in \sigma_p} \lambda_k e^{-t\lambda_k^2}
-\frac{\dsi}{4\pi }\int^{\infty}_{-\infty}\lambda e^{-t\lambda^2}
\tr(C_-(-i\lambda)C'_+(i\lambda))\ d\lambda\notag\\
& \quad \qquad \qquad \qquad +\frac{\dsi}{4\pi
}\int^{\infty}_{-\infty}\lambda e^{-t\lambda^2}\tr(C_+(-i\lambda)
C'_-(i\lambda))
\ d\lambda \notag \\
=&\quad \sum_{\lambda_k\in \sigma_p} \lambda_k e^{-t\lambda_k^2}
-\frac{\dsi}{4\pi }\int^{\infty}_{-\infty}\lambda e^{-t\lambda^2}
\tr_s(C_{\G}(\sigma_{+},-i\lambda) C'_{\G}(\sigma_{+},i\lambda))
\ d\lambda \notag \\
=&\quad  I_{\G}(K^o_t)+H_{\G}(K^o_t)+U_{\G}(K^o_t) \notag
\end{align}
where $d(\sigma_+)$ is the degree of $\sigma_+$ (although we know
$d(\sigma_+) = 2^{n-1}$, we will use the notation $d(\sigma_+)$
instead of $2^{n-1}$ since this indicates the origin of the
constant factor) and
$$
\tr_s(C_{\G}(\sigma_{+},-\nu)C'_{\G}(\sigma_{+},\nu))
:=\tr(C_-(-\nu)C'_+(\nu))-\tr(C_+(-\nu)C'_-(\nu)).
$$
For $I_{\G}(K_t)$, we have
$$
I_{\G}(K_t) =\text{Vol}(\G\bs G)\Bigr(\int^{\infty}_{-\infty}
\widehat{K}_t(\sigma_{+},i\lambda) p(\sigma_+,i\lambda)\ d\lambda
+\int^{\infty}_{-\infty} \widehat{K}_t(\sigma_-,i\lambda)
p(\sigma_-,i\lambda)\ d\lambda\Bigr)
$$
where $p(\sigma_{\pm},i\lambda)$ is the Plancherel measure which
is an even polynomial with respect to $\lambda$. The equalities
$\widehat{K}^e_t(\sigma_{\pm},i\lambda)=e^{-t\lambda^2}$ and
$p(\sigma_{+},i\lambda)=p(\sigma_{-},i\lambda)$ give
\begin{equation}
I_{\G}(K^e_t)=2\text{Vol}(\G\bs G) \int^{\infty}_{-\infty}
e^{-t\lambda^2}p(\sigma_{+},i\lambda)\ d\lambda
\end{equation}
and $\widehat{K}^o_t$ is odd with respect to $\lambda$ so that
\begin{equation}\label{e:oid}
I_{\G}(K^o_t)=0.
\end{equation}
It is well-known that
\begin{equation}\label{e:ezeta}
H_{\G}(K^e_t)=\frac{1}{\sqrt{4\pi t}}
\sum_{\gamma:\text{hyperbolic}}
l(C_{\gamma})j(\gamma)^{-1}D(\gamma)^{-1}\big(\,
\overline{\chi_{\sigma_{+}}(m_{\gamma})}+
\overline{\chi_{\sigma_{-}}(m_{\gamma})}\, \big)\,
e^{-\frac{l(C_{\gamma})^2}{4t}},
\end{equation}
\begin{equation}\label{e:ozeta}
H_{\G}(K^o_t)=\frac{2\pi i}{{(4\pi t)}^{\frac 32}}
\sum_{\gamma:\text{hyperbolic}}
l(C_{\gamma})^2j(\gamma)^{-1}D(\gamma)^{-1}\big(\,
\overline{\chi_{\sigma_{+}}(m_{\gamma})}-
\overline{\chi_{\sigma_{-}}(m_{\gamma})}\, \big)\,
e^{-\frac{l(C_{\gamma})^2}{4t}}
\end{equation}
where $l(C_{\gamma})$ is the length of the closed geodesic
$C_{\gamma}$, $j(\gamma)$ is the positive integer such that
$\gamma=\gamma_0^{j(\gamma)}$ for a primitive element $\gamma_0$,
$D(\gamma)=e^{n
l(C_{\gamma})}\vert\det(\text{Ad}(a_{\gamma}m_{\gamma})^{-1} -
I\vert_{\frak{n}})\vert$ for the element $a_{\gamma}m_{\gamma}\in
A^{+}M$ which is conjugate to $\gamma$ and $\chi_{\sigma}$ is the
character of $\sigma$.

\section{Unipotent terms}\label{s3}

In this section we compute unipotent terms $U_{\G}({K}^e_t)$,
$U_{\G}({K}^o_t)$. We employ the formula obtained by Hoffmann (see
\cite {H}) to compute these terms explicitly.  By this explicit
computation and Proposition \ref{p:frie}, it follows that
$U_{\Gamma}(K^o_t)=0$. This simplifies many steps related to the
application of the Selberg trace formula for $K^o_t$.

For a real rank $1$ group $G$, the unipotent term for a right
$K$-finite function $f$ in $\mathcal{C}^p(G)$, $0<p<1$  is given
by
\begin{equation}\label{e:U}
U_{\G}(f)=\frac{1}{\vert\alpha_1\vert}\{C_1(\G)T_1(f)+C_2(\G)T_2(f)+C'_1(\G)T'_1(f)\},
\end{equation}
(see theorem in p.299 of \cite{OW}). Here $\alpha_1$ is the unique
simple positive root for $(\frak{g}, \frak{a})$, the constants
$C_1(\G), C_2(\G), C'_1(\G)$ which depend on $\G$ are computed in
\cite{D} and
\begin{align*}
T_1(f)&=\frac{1}{A(\frak{n}_1)}\int_{N}\int_{K} f(k^{-1}nk)\, dkdn\\
T_2(f)&=\frac{\vert\alpha_1\vert}{2}\big\{\int_{G/G_{n_0}}
f(xn_0x^{-1})\, dx
+\int_{G/G_{n_0^{-1}}} f(xn_0^{-1}x^{-1})\, dx \big\} \\
T'_1(f)&=\frac{m_1+2m_2}{A(\frak{n}_1)}\int_{N_1}\int_{N_2}\int_K
f(k^{-1}n_1n_2k)\log \vert \log(n_1)\vert\, dk dn_1 dn_2
\end{align*}
where $n_0$ is a representative of the nontrivial unipotent orbit
in $\frak{n}_2$, $m_i=\text{dim} (\frak{n_i})$ and $A(\frak{n}_1)$
is the volume of the unit sphere in $\frak{n}_1$ and $G_n$ is the
centralizer of an element $n$ in $G$.

In our case, $G=\Spin(2n+1,1)$, the second term $T_2(f)=0$ and
$m_2=0$ in the third term $T'_1(f)$ since $N=N_1$. The first term
$T_1$ for $K_t$ is given by
$$
T_1(K_t)=\frac{1}{2\pi A(\frak{n}_1)}
\Bigr(\int^{\infty}_{-\infty} \widehat{K}_t(\sigma_+,i\lambda)\
d\lambda +\int^{\infty}_{-\infty}
\widehat{K}_t(\sigma_{-},i\lambda)\ d\lambda\Bigr).
$$
We have $\widehat{K}^e_t(\sigma_{\pm},i\lambda)=e^{-t\lambda^2}$
and $\widehat{K}^o_t(\sigma_{\pm},i\lambda)=\pm\lambda
e^{-t\lambda^2}$, which implies
\begin{equation}\label{e:unip0}
T_1(K^e_t)=\frac{1}{\pi A(\frak{n}_1)} \int^{\infty}_{-\infty}
e^{-t\lambda^2} \ d\lambda \quad\text{and}\quad T_1(K^o_t)=0.
\end{equation}

The third term $T'_1(f)$ is more complicated and we need to
introduce some notations. Let $T_M$ be a Cartan subgroup in $M$,
so that $T=A\cdot T_M$ is a Cartan subgroup of $G$. Let $\Sigma_M$
denote the set of positive roots for
$(\frak{m}_{\mathbb{C}},\frak{t_m}_{\mathbb{C}})$. Let
$\rho_{\Sigma_M}$ be the half sum of elements in $\Sigma_M$. We
denote by $\Sigma_A$ the set of positive roots of
$(\frak{g}_{\mathbb{C}}, \frak{t}_{\mathbb{C}})$ which do not
vanish on $\frak{a}_{\mathbb{C}}$. The union of $\Sigma_M$ with
$\Sigma_A$ gives the set of positive roots for
$(\frak{g}_{\mathbb{C}}, \frak{t}_{\mathbb{C}})$ denoted by
$\Sigma_G$. Let $H_{\alpha}\in \frak{t}_{\mathbb{C}}$ be the
coroot corresponding to $\alpha\in\pm\Sigma_G$, that is,
$\alpha(H_{\alpha})=2,\alpha'(H_{\alpha})\in\mathbb{Z}$ for all
$\alpha,\alpha'\in\pm\Sigma_G$ and $\Pi:=\prod_{\alpha\in\Sigma_M}
H_{\alpha} $ is an element of the symmetric algebra
$S(\frak{t_m}_{\mathbb{C}})$. We denote the simple reflection
corresponding to $\alpha\in \Sigma_G$ by $s_{\alpha}$. Following
the corollary on p.96 of \cite{H}, we put
$$
I(f):=\frac{1}{2\pi}\sum_{\sigma\in\hat{M}}\int_{-\infty}^{\infty}
\Omega(\sigma,-i\lambda)\hat{f}(\sigma,i\lambda)\ d\lambda
$$
where
\begin{equation*}
\Omega(\sigma,i\lambda)=2 d(\sigma)\psi(1)
-\frac{1}{2}\sum_{\alpha\in\Sigma_A}
\lambda(H_{\alpha})\frac{\Pi(s_{\alpha}\lambda_{\sigma})}{\Pi(\rho_{\Sigma_M})}
\times\Bigr(\psi(1+\lambda_{\sigma}(H_{\alpha}))+\psi(1-\lambda_{\sigma}(H_{\alpha}))\Bigr)
\end{equation*}
where $d(\sigma)$ is the degree of $\sigma$, $\psi$ is the
logarithmic derivative of the Gamma function and
$\lambda_{\sigma}-\rho_{\Sigma_M}$ is the highest weight of
$(\sigma,i\lambda)\in \hat{M}\times i\frak{a}$.

The principal series representation
$(\pi_{\sigma,\nu},\mathcal{H}_{\sigma,\nu})$ depends on the
parabolic subgroup $P$ and we denote its dependence on $P$ by
$(\pi_{\sigma,\nu}(P),\mathcal{H}_{\sigma,\nu}(P))$. The
intertwining operator
$$
J_{\bar{P}\vert P}(\sigma,\nu):\mathcal{H}_{\sigma,\nu}(P)\to
\mathcal{H}_{\sigma,\nu}(\bar{P})
$$
is defined by
$$
(J_{\bar{P}\vert P}(\sigma,\nu)\phi)=\int_{\bar{N}}\phi(x\bar{n})\
d\bar{n}
$$
and satisfies
$$
J_{\bar{P}\vert P}(\sigma,\nu) \pi_{\sigma,\nu}(P)
=\pi_{\sigma,\nu}(\bar{P}) J_{\bar{P}\vert P}(\sigma,\nu).
$$
The restriction to $K$ defines an isomorphism from
$\mathcal{H}_{\sigma,\nu}(P)$ to $\mathcal{H}_{\sigma}(P)$. Then
$J_{\bar{P}\vert P}(\sigma,\nu)$ can be considered as a family of
operators from $\mathcal{H}_{\sigma}(P)$ to
$\mathcal{H}_{\sigma}(\bar{P})$. Here $\mathcal{H}_{\sigma}(P)$ is
the space of all measurable functions $v:K \to {H}_{\sigma}$ such
that
$$
v(km)=\sigma(m)^{-1}v(k)
$$
for all $m\in M$, $k\in K$. Let
$$
J_P(\sigma,\nu:f):=-\Tr(\pi_{\sigma,\nu}(f)J_{\bar{P}\vert
P}(\sigma,\nu)^{-1}\partial_{\nu} J_{\bar{P}\vert P}(\sigma,\nu))
$$
where the derivative $\partial_{\nu}$ is taken with respect to
$\nu$ for the family of operators $J_{\bar{P}\vert P}(\sigma,\nu)$
acting on $\mathcal{H}_{\sigma}(P)$. Then there exists the Harish
Chandra $C$-function $C_{\tau}(\sigma,\nu)$ such that
$$
T_{\tau} J_{\bar{P}\vert P}(\sigma,\nu)^{-1}\partial_{\nu}
J_{\bar{P}\vert P}(\sigma,\nu) \ = \
C_{\tau}(\sigma,\nu)^{-1}\partial_{\nu}C_{\tau}(\sigma,\nu)
T_{\tau}
$$
where $T_{\tau}$ is the projection to the $\tau$-isotypic
component of $\mathcal{H}_{\sigma,\nu}(P)$. We refer to \cite{EKM}
for more detail. We have the following proposition for
$C_{\tau_n}(\sigma_{\pm},\nu)$,

\begin{prop}\label{p:unip1} For the half spin representation
$\sigma_{\pm}$ of $\Spin(2n)$,
$$
C_{\tau_n}(\sigma_{+}, i\lambda)= C_{\tau_n}(\sigma_{-},i\lambda)
=  \frac{(2n-1)!}{(n-1)!}\frac{\Gamma(i\lambda +\frac 12)}
{\Gamma(i\lambda +n +\frac 12)}\mathrm{Id}.
$$
\end{prop}

\begin{proof}
This follows from theorem 8.2 in \cite{EKM}.
\end{proof}

By the equality (8) (see also (49)) in \cite{H}, the weighted
orbital integral $T'_1 $ is given by
$$
T'_1(f)=
\frac{m_1}{A(\frak{n}_1)}\Bigr(I(f)+\frac{1}{2\pi}\text{p.v.}\sum_{\sigma\in
\hat{M}}d(\sigma) \int^{\infty}_{-\infty} J_P(\sigma, i\lambda
:f)\ d\lambda
+\sum_{\sigma\in\hat{M}}d(\sigma)\frac{n(\sigma)}{2}\hat{f}(\sigma,0)\Bigr)
$$
where p.v.\ means the Cauchy principal value and $2n(\sigma)$ is
the order of the zero of $p(\sigma,\nu)$ at $\nu=0$. Using that
$n(\sigma_{\pm})=0$ and Proposition \ref{p:unip1} we obtain
$$
T'_1(K^e_t)=\frac{m_1}{2\pi A(\frak{n}_1)}
\int^{\infty}_{-\infty}e^{-t\lambda^2}\Bigr(\Omega(\sigma_+,-i\lambda)+\Omega(\sigma_-,-i\lambda)
-2d(\sigma_+)\partial_{i\lambda} \log
C_{\tau_n}(\sigma_+,i\lambda)\Bigr)\ d\lambda,
$$
$$
T'_1(K^o_t)=\frac{m_1}{2\pi A(\frak{n}_1)}
\int^{\infty}_{-\infty}\lambda e^{-t\lambda^2}\Bigr
(\Omega(\sigma_+,-i\lambda)- \Omega(\sigma_-,-i\lambda)\Bigr)\
d\lambda .
$$
Now we consider $\Omega(\sigma_{\pm},i\lambda)$ in the following
proposition.

\begin{prop}\label{p:unip2}
For the half spin representation $\sigma_{\pm}$ of $\Spin(2n)$, we
have
\begin{multline*}
\Omega(\sigma_{+},i\lambda)\ =\ \Omega(\sigma_{-},i\lambda)\\
 =\
-\frac{d(\sigma_{+})}{2}\Bigr(\psi(i\lambda-n+\frac 12)
+\psi(-i\lambda-n+\frac 12) +\psi(i\lambda+ \frac 12)
+\psi(-i\lambda +\frac 12)\Bigr) +P^n(\lambda)
\end{multline*}
where $P^n(\lambda)$ is an even polynomial of degree $(2n-4)$ for
$n\geq 2$ and $P^1(\lambda)$ is a constant.
\end{prop}

\begin{proof}
The $n=1$ case can be computed in the same way as for the $n\geq
2$ cases, hence we may assume that $n\geq 2$ in the following
proof. The highest weight of the half spin representation
$\sigma_{\pm}$ of $\Spin(2n)\subset \Spin(2n+1)$ is given by
$$
\frac{1}{2}(e_2+e_3+\cdots +e_n \pm e_{n+1})
$$
with respect to the standard basis $\{e_i\}$. This implies that
$$
i\lambda e_1 + \sigma_{\pm} + \rho_{\Sigma_M} =i\lambda e_1 +
(n-\frac{1}{2}) e_2 +(n-\frac{3}{2}) e_3 + \cdots + \frac{3}{2}
e_n \pm \frac{1}{2} e_{n+1}.
$$
The positive roots $\alpha\in\Sigma_A$ are given by $e_1 - e_j,
e_1 + e_j$ for $2\leq j \leq n+1$. Then we can see that
$s_{(e_1-e_j)}( i\lambda e_1 + \sigma_{\pm} + \rho_{\Sigma_M}) $
is
\begin{align*}
& i\lambda e_j+ (n-\frac{1}{2}) e_2 +\cdots
+ (n-j+\frac{3}{2})e_1+\cdots+\frac{3}{2} e_n \pm \frac{1}{2} e_{n+1} \quad \text{if}\quad 2\leq j \leq n, \\
&i\lambda e_{n+1}+ (n-\frac{1}{2}) e_2 +\cdots +\frac{3}{2} e_n
\pm \frac{1}{2} e_{1} \quad \text{if}\quad j=n+1
\end{align*}
and $s_{(e_1+ e_j)}( i\lambda e_1 + \sigma_{\pm} +
\rho_{\Sigma_M}) $ is
\begin{align*}
&-i\lambda e_j+ (n-\frac{1}{2}) e_2 +\cdots
- (n-j+\frac{3}{2})e_1+\cdots+\frac{3}{2} e_n \pm \frac{1}{2} e_{n+1} \quad \text{if}\quad 2\leq j \leq n,\\
&-i\lambda e_{n+1}+ (n-\frac{1}{2}) e_2 +\cdots +\frac{3}{2} e_n
\mp \frac{1}{2} e_{1} \quad\text{if}\quad j=n+1.
\end{align*}
These give us the formula
\begin{align}\label{p:poly}
&\Pi(s_{(e_1\pm e_j)}( i\lambda e_1 + \sigma_{\pm} + \rho_{\Sigma_M}))\notag \\
=&
C_j(\lambda^2+(n-\frac 12)^2)\cdot(\lambda^2+(n-\frac 32)^2)\cdots(\lambda^2+(n-j+\frac 52)^2)\\
&\cdot (-\lambda^2-(n-j+\frac 12)^2)\cdots(-\lambda^2-(\frac
32)^2)\cdot(-\lambda^2-(\frac 12)^2)\notag
\end{align}
where
\begin{align*}
C_j&=\prod_{ 1\leq k<l\leq n, k\neq j,l\neq j}
((n-k+\frac{1}{2})^2-(n-l+\frac{1}{2})^2) \\
&=\prod_{ 1\leq k<l\leq n, k\neq j,l\neq j}(l-k)(2n-k-l+1)
\end{align*}
for $2 \leq j \leq n+1$. In particular, we obtain
\begin{align*}
\Pi(s_{e_1\pm e_j}(i\lambda e_1 +\sigma_{+} + \rho_{\Sigma_M}))=
\Pi(s_{e_1\pm e_j}(i\lambda e_1 +\sigma_{-} + \rho_{\Sigma_M}))\\
\Pi(s_{e_1-e_j}(i\lambda e_1 +\sigma_{\pm} + \rho_{\Sigma_M}))
=\Pi(s_{e_1+e_j}(i\lambda e_1 +\sigma_{\pm} + \rho_{\Sigma_M}))
\end{align*}
for $2\leq j\leq n+1$. From now on, we denote by $P_j(\lambda)$
the polynomial of $\lambda$ in (\ref{p:poly}) for $2\leq j\leq
n+1$. Note that $P_j(\lambda)$ is an even polynomial of degree
$2(n-1)$. On the other hand, $(i\lambda e_1 +\sigma_{\pm} +
\rho_{\Sigma_M})(H_{\alpha})$ is given by
\begin{align*}
&i\lambda -(n-j+\frac{3}{2}) \quad\text{if}\quad \alpha= e_1-
e_{j},
\quad 2\leq j\leq n \\
&i\lambda \mp \frac{1}{2} \quad\text{if}\quad \alpha=e_1-e_{n+1} \\
&i\lambda +(n-j+\frac{3}{2}) \quad\text{if}\quad \alpha=
e_1+e_{j},
 \quad 2\leq j\leq n \\
&i\lambda \pm \frac{1}{2} \quad\text{if}\quad \alpha=e_1+e_{n+1}.
\end{align*}
Then the pair $\Bigr(\psi(1+(i\lambda e_1 +\sigma_{\pm} +
\rho_{\Sigma_M})(H_{\alpha})), \psi(1-(i\lambda e_1 +\sigma_{\pm}
+ \rho_{\Sigma_M})(H_{\alpha}))\Bigr)$ is given by
\begin{align*}
&\psi(i\lambda-n+j-\frac 12),\ \ \psi(-i\lambda+n-j+\frac 52)
\quad\text{for}\quad  e_1-e_j,
\quad 2\leq j\leq n\\
&\psi(i\lambda +\frac 12),\ \ \psi(-i\lambda+\frac 32) \quad \text{for}\quad  e_1-e_{n+1}, \quad \sigma=\sigma_{+}\\
&\psi(i\lambda+\frac 32),\ \ \psi(-i\lambda +\frac 12) \quad \text{for} \quad e_1-e_{n+1}, \quad \sigma=\sigma_-\\
&\psi(i\lambda+n-j+\frac 52),\ \ \psi(-i\lambda-n+j-\frac 12)
\quad\text{for}\quad e_1+e_{j},
\quad  2\leq j\leq n \\
&\psi(i\lambda+\frac 32), \ \ \psi(-i\lambda +\frac 12) \quad \text{for}\quad e_1+e_{n+1}, \quad \sigma=\sigma_+\\
&\psi(i\lambda+\frac 12), \ \ \psi(-i\lambda +\frac 32) \quad
\text{for}\quad e_1+e_{n+1}, \quad \sigma=\sigma_-.
\end{align*}
Comparing the two sets
\begin{align*}
&\Bigr\{\Pi(s_{\alpha}(i\lambda e_1 +\sigma_{+} +
\rho_{\Sigma_M})), 1+(i\lambda e_1 +\sigma_{+} +
\rho_{\Sigma_M})(H_{\alpha}),
1-(i\lambda e_1 +\sigma_{+} + \rho_{\Sigma_M})(H_{\alpha})\Bigr\},\\
&\Bigr\{\Pi(s_{\alpha}(i\lambda e_1 +\sigma_{-} +
\rho_{\Sigma_M})), 1+(i\lambda e_1 +\sigma_{-} +
\rho_{\Sigma_M})(H_{\alpha}), 1-(i\lambda e_1 +\sigma_{-} +
\rho_{\Sigma_M})(H_{\alpha})\Bigr\}
\end{align*}
we see that they are equal to each other, so that
$\Omega(\sigma_{+},i\lambda)=\Omega(\sigma_{-},i\lambda)$. We now
compute the exact form of
$\Omega(\sigma_{+},i\lambda)=\Omega(\sigma_{-},i\lambda)$. Using
the relations $ \psi(z+1)= \frac 1 z +\psi(z)$, we have
\begin{align*}
&\psi(i\lambda -n+j-\frac 12) + \psi(-i\lambda -n+j-\frac 12) +
\psi(i\lambda+n-j+\frac 52)+ \psi(-i\lambda+n-j+\frac 52) \\
=&\frac 1{\lambda^2+(\frac12)^2}+\cdots +\frac {2(n-j+\frac
12)}{\lambda^2+(n-j+\frac 12)^2}+\frac {-2(n-j+\frac
52)}{\lambda^2+(n-j+\frac 52)^2}+\cdots + \frac{-2(n-\frac
12)}{\lambda^2+(n-\frac 12)^2}\\
&+\psi(i\lambda-n+\frac 12) +\psi(-i\lambda-n+\frac 12)
+\psi(i\lambda+\frac 12) +\psi(-i\lambda+\frac 12) .
\end{align*}
Now using the formula (see the last line of p.95 in \cite{H}),
$$\sum_{\alpha\in \Sigma_A}
\Pi(s_{\alpha}\lambda_{\sigma})=2\Pi(\lambda_{\sigma}) , $$ we
decompose
$$
\frac{1}{2}\sum_{\alpha\in\Sigma_A}
\frac{\Pi(s_{\alpha}\lambda_{\sigma})}{\Pi(\rho_{\Sigma_M})}
\times\Bigr(\psi(1+\lambda_{\sigma}(H_{\alpha}))+\psi(1-\lambda_{\sigma}(H_{\alpha}))\Bigr)
$$
into
$$
\frac{d(\sigma_{\pm})}{2}\Bigr(\psi(i\lambda-n+\frac 12)
+\psi(-i\lambda-n+\frac 12) +\psi(i\lambda+ \frac 12)
+\psi(-i\lambda +\frac 12)\Bigr)
$$
and
$$
\frac{1}{4\Pi(\rho_{\Sigma_M})}\sum_{j=2}^{n+1}P_j(\lambda)
R_j(\lambda)
$$
where
$$
R_j(\lambda)=\frac 1{\lambda^2+(\frac12)^2}+\cdots +\frac
{2(n-j+\frac 12)}{\lambda^2+(n-j+\frac 12)^2}+\frac {-2(n-j+\frac
52)}{\lambda^2+(n-j+\frac 52)^2}+\cdots + \frac{-2(n-\frac
12)}{\lambda^2+(n-\frac 12)^2}.
$$
We can see, from the definitions of $P_j(\lambda)$ and
$R_j(\lambda)$, that $P_j(\lambda)R_j(\lambda)$ is an even
polynomial of degree $2n-4$. This ends the proof.
\end{proof}

For the constants in \eqref{e:U}, comparing proposition 6.2 in
\cite{BM} with  theorem 2 in \cite{D}, we obtain
\begin{equation}\label{e:counip}
\frac{C'_1(\G)}{\vert\alpha_1\vert}\frac{m_1}{A(\frak{n}_1)}=\frac{\kappa}{2}
.
\end{equation} Now we have the following corollary.

\begin{cor}\label{c:ounip}
$$
U_{\G}(K^o_t)=0, \qquad
U_{\G}(K^e_t)=\frac{2}{2\pi}\int^{\infty}_{-\infty}
e^{-t\lambda^2} (P_U(\lambda)+Q(\lambda))\ d\lambda
$$
where $P_U(\lambda)$ is an even polynomial of degree $(2n-4)$ and
$$ Q(\lambda)= -\kappa \frac{d(\sigma_{+})}{2}\Bigr(\psi(i\lambda+
\frac 12) +\psi(-i\lambda +\frac 12) \Bigr).$$

\end{cor}

\begin{proof}
The first claim follows easily from (\ref{e:unip0}),
(\ref{e:counip}), Proposition \ref{p:unip1} and \ref{p:unip2}. For
the second claim, $Q(\lambda)$ is a priori given by
\begin{align*}
&-\frac{\kappa d(\sigma_{+})}{2}\Bigr(\, \frac 12 \big(\,
\psi(i\lambda-n+\frac 12)
+\psi(-i\lambda-n+\frac 12)\, \big)\\
&\qquad \qquad - \frac 12 \big(\, \psi(i\lambda+n+\frac 12)
+\psi(-i\lambda+n+\frac 12)\, \big)+ \big(\, \psi(i\lambda+ \frac
12) +\psi(-i\lambda +\frac 12)\, \big)\, \Bigr).
\end{align*}
If we use the relation $\psi(z+1)=\frac 1z +\psi(z)$, then we can
reduce the above formula to the claimed one for $Q(\lambda)$.
\end{proof}

\section{Relative Traces and Spectral Sides}\label{s4}

In this section we study relations of the relative traces with the
spectral sides of the Selberg trace formulas applied to the test
functions $K^e_t,K^o_t$. A formula of this type was proved by
M\"uller for the similar cases  in \cite{Mu1}, \cite{Mu3}.
Following \cite{Mu1}, \cite{Mu3}, we prove the corresponding
formula for Dirac operators acting on spinor bundles over
hyperbolic manifolds with cusps.

First, let us observe that $L^2(X,E)$ can be identified with the
space
$$
(L^2(\G\bs G)\otimes V_{\tau_{n}})^K=\{ \ f\in L^2(\G\bs G)\otimes
V_{\tau_{n}}\ |\ f(xk)=\tau_{n}(k)^{-1}f(x) \quad \text{for}\quad
k\in K, x\in G \ \}.
$$
The decomposition of $L^2(\G\bs G)=L^2_d(\G\bs G)\oplus
L^2_c(\G\bs G)$ allows us to decompose
\begin{equation}\label{e:spec-decom}
(L^2(\G\bs G)\otimes V_{\tau_{n}})^K = (L_d^2(\G\bs G)\otimes
V_{\tau_{n}})^K \oplus (L_c^2(\G\bs G) \otimes V_{\tau_{n}})^K.
\end{equation}
The continuous part $ ( L_c^2(\G\bs G) \otimes V_{\tau_{n}})^K$ is
spanned by the wave packets with the Eisenstein series, so we need
to know how $D$ acts on the Eisenstein series. For $\Phi \in (
\mathcal{H}_{\Gamma}(\sigma_+,\nu)\otimes V_{\tau_n})^K$, the
Eisenstein series attached to $\Phi$ is defined by
$$
E(\Phi : \nu :x):= \sum_{\gamma\in \Gamma\cap P\bs \Gamma} \Phi
(\gamma x),
$$
which is defined a priori for $\Re(\nu)\gg 0$ and has the
meromorphic extension over $\mathbb{C}$. Assume that
$\Phi_{j,\pm}$ is in the $\pm i$-eigenspace of $c(H)$ on $(
\mathcal{H}_{\Gamma}(\sigma_+,\nu)\otimes V_{\tau_n})^K$ for
$j=1,\cdots,d(\sigma_+)$. The Eisenstein series $E(\Phi_{j,\pm}:
i\lambda :x)$ for $\nu=i\lambda$ satisfies
\begin{equation}\label{e:DdE}
D E(\Phi_{j,\pm}:i\lambda: x)=\pm \lambda
E(\Phi_{j,\pm}:i\lambda:x).
\end{equation}
For $\phi\in C^{\infty}_0(X,E)$, the decomposition
(\ref{e:spec-decom}) provides us with the formula
\begin{align}\label{e:spec-decom1}
\phi(x)&=\sum_k(\phi,\phi_k)\phi_k(x)\\
&\quad +\frac{1}{4\pi}\sum_{j=1}^{d(\sigma_+)}
\int^{\infty}_{-\infty} E( \Phi_{j,+}:i\lambda:x)\int_{X}
E(\Phi_{j,+}:-i\lambda:y)\phi(y)\ dy
\ d\lambda \notag \\
&\quad
+\frac{1}{4\pi}\sum_{j=1}^{d(\sigma_+)}\int^{\infty}_{-\infty} E(
\Phi_{j,-}:i\lambda:x)\int_{X} E(\Phi_{j,-}:-i\lambda:y)\phi(y)\
dy \ d\lambda \notag
\end{align}
where  $\{\phi_k\}$ is an orthonormal basis of $(L_d^2(\G\bs
G)\otimes V_{\tau_{n}})^K $. Since $\Dd$ preserves the
decomposition of (\ref{e:spec-decom}), we may assume that each
$\phi_k$ is an eigensection of $\Dd$. Therefore (\ref{e:DdE}) and
(\ref{e:spec-decom1}) imply
\begin{align*}
\Dd\phi(x)&=
\sum_{\lambda_k} \lambda_k (\phi,\phi_k) \phi_k(x)\\
&\quad
+\frac{1}{4\pi}\sum_{j=1}^{d(\sigma_+)}\int^{\infty}_{-\infty}
\lambda E( \Phi_{j,+}:i\lambda:x)
\int_X E( \Phi_{j,+}:-i\lambda:y)\phi(y)\ dy \ d\lambda\\
&\quad
-\frac{1}{4\pi}\sum_{j=1}^{d(\sigma_+)}\int^{\infty}_{-\infty}
\lambda E( \Phi_{j,-}:i\lambda:x) \int_X E(
\Phi_{j,-}:-i\lambda:y)\phi(y)\ dy \ d\lambda .
\end{align*}
From this we can see that the following equalities hold in the
distributional sense:
\begin{align}\label{e:trace-e}
e^{-t\Dd^2}(x:y)&=\sum_{\lambda_k\in \sigma_p} e^{-t\lambda^2_k}\phi_k(x)\otimes \phi^*_k(y)\\
&\quad
+\frac{1}{4\pi}\sum_{j=1}^{d(\sigma_+)}\int^{\infty}_{-\infty}
e^{-t\lambda^2} E( \Phi_{j,+}:i\lambda:x)
\otimes E( \Phi_{j,+}:-i\lambda:y)^{*} \ d\lambda \notag \\
&\quad
+\frac{1}{4\pi}\sum_{j=1}^{d(\sigma_+)}\int^{\infty}_{-\infty}
e^{-t\lambda^2} E( \Phi_{j,-}:i\lambda:x) \otimes E(
\Phi_{j,-}:-i\lambda:y)^{*} \ d\lambda, \notag
\end{align}
\begin{align}\label{e:trace-o}
\Dd e^{-t\Dd^2}(x:y)&=\sum_{\lambda_k\in\sigma_p} \lambda_k e^{-t\lambda^2_k}\phi_k(x)\otimes \phi^*_k(y)\\
&\quad
+\frac{1}{4\pi}\sum_{j=1}^{d(\sigma_+)}\int^{\infty}_{-\infty}\lambda
e^{-t\lambda^2} E( \Phi_{j,+}:i\lambda:x)
\otimes E( \Phi_{j,+}:-i\lambda:y)^{*}\ d\lambda \notag \\
&\quad
-\frac{1}{4\pi}\sum_{j=1}^{d(\sigma_+)}\int^{\infty}_{-\infty}
\lambda e^{-t\lambda^2} E( \Phi_{j,-}:i\lambda:x) \otimes E(
\Phi_{j,-}:-i\lambda:y)^{*}\ d\lambda . \notag
\end{align}

Recalling the decomposition \eqref{decompX} of $X$, we consider
the subset $W_{i,R}:=[R,\infty)\times N_i$ in $W_i$. We choose
$R_0$ such that $W_{i,R}$'s are disjoint each other for $R\ge
R_0$. We put $W_R:=\cup^{\kappa}_{i=1} W_{i,R}$ and $X_R:= X
-W_{R}$ for $R\ge R_0$. The constant term of the Eisenstein series
$E(\Phi:i\lambda:x)$ over $W_{i}$ has the form
\begin{equation}\label{e:const}
e^{(-i\lambda+n)r}\Phi_i + e^{(i\lambda+n)r}\big( (
C_{\Gamma}(\sigma_{+},i\lambda)\otimes {\Id}) \Phi\big)_i.
\end{equation}
Here $\Phi_i$ is the component of $\Phi$ over $W_{i}$ and note
that the operator $C_{\Gamma}(\sigma_{+},i\lambda)\otimes \Id$
acts on $( \mathcal{H}_{\Gamma}(\sigma_+,i\lambda)\otimes
V_{\tau_{n}})^K$. From now on, we assume that $||\Phi_i||=1$ for
$i=1,\cdots, \kappa$. We now discuss the {\it Maass-Selberg
relation} in our context.

\begin{prop}\label{p:maass}{\rm(Maass-Selberg)} We have
\begin{align}\label{e:Eisen}
\int_{X_R} \vert E(\Phi_{j,\pm}:i\lambda:x)\vert^2 dx \ = &\
2\kappa R - \tr\big(\,C_{\pm}(i\lambda) C'_{\mp}(-i\lambda)\,\big)
+O(e^{-cR}) \\
=&\ 2\kappa R - \tr\big(\,C_{\mp}(-i\lambda) C'_{\pm}(i\lambda)\,
\big) +O(e^{-cR}) \notag
\end{align}
where $c$ is a positive constant.
\end{prop}

\begin{proof}
We will consider only the case of $\Phi_+$ so we use the notation
$\Phi$ instead of $\Phi_+$ in the following proof. The case of
$\Phi_-$ can be done in the same way. It follows from Green's
formula that
\begin{align}\label{e:maass1}
&(\lambda-\lambda')\langle\, E(\Phi:i\lambda:x),
E(\Phi:i\lambda':x)
\, \rangle_{X_R} \notag \\
=&\ \langle\, D E(\Phi:i\lambda:x), E(\Phi:i\lambda':x) \,
\rangle_{X_R}-\langle\, E(\Phi:i\lambda:x), D E(\Phi:i\lambda':x)
\, \rangle_{X_R}  \\
=&\ \langle\, c(H) E(\Phi:i\lambda:x), E(\Phi:i\lambda':x) \,
\rangle_{\partial X_R}. \notag
\end{align}
By (\ref{e:const}), we have
\begin{align*}
&\ \langle\, c(H) E(\Phi:i\lambda:x), E(\Phi:i\lambda':x)
\, \rangle_{\partial W_{i,R}} \notag \\
=&\ i\, \langle\, \ e^{(-i\lambda+n)R}\Phi_{i}\ ,\
e^{(-i\lambda'+n)R}\Phi_{i}\, \rangle_{\partial
W_{i,R}} \\
& - i\, \langle \, e^{(i\lambda+n)R}(
(C_{\Gamma}(\sigma_{+},i\lambda)\otimes \text{Id})\Phi)_i \ , \
e^{(i\lambda'+n)R}( (C_{\Gamma}(\sigma_{+},i\lambda')\otimes
\text{Id})\Phi )_i\, \rangle_{\partial
W_{i,R}} + O(e^{-cR}) \notag \\
=&\  i\, e^{-i(\lambda-\lambda')R} - i\,
e^{i(\lambda-\lambda')R}\sum^{\kappa}_{k=1}
C_{+}(i\lambda)_{ik}C_{-}(-i\lambda')_{ki} +O(e^{-cR}) \notag
\end{align*}
where $C_{\pm}(i\lambda)_{ik}$ is a component of
$C_{\pm}(i\lambda)$. Now we use the functional equation
$$
C_{\G}(\sigma_+,i\lambda)C_{\G}(\sigma_+,-i\lambda)=\left (
\begin{matrix} 0 & C_-(i\lambda) \\
C_+(i\lambda) & 0
\end{matrix}\right )
\left (
\begin{matrix} 0 & C_-(-i\lambda) \\
C_+(-i\lambda) & 0
\end{matrix}\right )=\text{Id} ,$$
which implies the equality
\begin{align}\label{e:maass2}
&\  \langle\, c(H) E(\Phi:i\lambda:x), E(\Phi:i\lambda':x)
\, \rangle_{\partial W_{i,R}} \notag \\
= &\ i\, e^{-i(\lambda-\lambda')R} - i\, e^{i(\lambda-\lambda')R}  \\
& -  i\, e^{i(\lambda-\lambda')R}\sum^{\kappa}_{k=1}
C_{+}(i\lambda)_{ik}C_{-}(-i\lambda')_{ki}  \notag \\
& + i\, e^{i(\lambda-\lambda')R}\sum^{\kappa}_{k=1}
C_{+}(i\lambda)_{ik}C_{-}(-i\lambda)_{ki} +O(e^{-cR}).\notag
\end{align}
Combining (\ref{e:maass1}) and (\ref{e:maass2}), we have
\begin{align}\label{e:mss}
& (\lambda-\lambda') \langle\, E(\Phi:i\lambda:x),
 E(\Phi:i\lambda':x) \, \rangle_{X_R} \ = \ i\, \kappa\,
e^{-i(\lambda-\lambda')R} - i\, \kappa\, e^{i(\lambda-\lambda')R} \\
& \qquad  +i\, e^{i(\lambda-\lambda')R}\sum^{\kappa}_{i,k=1}
C_{+}(i\lambda)_{ik}\big(C_{-}(-i\lambda)_{ki}
-C_{-}(-i\lambda')_{ki}\big) +O(e^{-cR}) . \notag
\end{align}
If we pass to the limit $\lambda\to \lambda'$ after dividing each
side of (\ref{e:mss}) by $\lambda-\lambda'$, we get the equality
(\ref{e:Eisen}).

\end{proof}

The next task is to consider the corresponding formula for
$\Dd_0(i)$. We fix a Lagrangian subspace $L$ in $V_{\tau_n}$ as
before. For an orthonormal basis $\{ \psi_{j}\} \subset L$, we
define
$$
\phi_{j,+}:=\frac{1}{\sqrt{2}}\biggl(\psi_{j}-ic(H)\psi_{j}\biggr)
, \quad
\phi_{j,-}:=\frac{1}{\sqrt{2}}\biggl(\psi_{j}+ic(H)\psi_{j}\biggr)
$$
and
$$
e(\phi_{j,+}:i\lambda:x):=
e^{(n-i\lambda)r}\phi_{j,+}+e^{(n+i\lambda)r}\phi_{j,-},
$$
$$
e(\phi_{j,-}:i\lambda:x):=
e^{(n-i\lambda)r}\phi_{j,-}+e^{(n+i\lambda)r}\phi_{j,+} .
$$
Note that  $e(\phi_{j,\pm}:i\lambda:x)$ lies in
$C^{\infty}(\mathbb{R}^+, V_{\tau_n}, L ):=\{\ \phi\in
C^{\infty}(\mathbb{R}^{+},E_0)\ | \ \phi(0)\in L\ \}$ and
\begin{align*}
& c(H)(\frac{d}{dr}-n\,\Id)\, e(\phi_{j,+}:i\lambda:x)=\lambda\, e(\phi_{j,+}:i\lambda:x), \\
& c(H)(\frac{d}{dr}-n\, \Id)\, e(\phi_{j,-}:i\lambda:x)=-\lambda\,
e(\phi_{j,-}:i\lambda:x).
\end{align*}
As when we introduced $\Dd_0(i)$,  we can regard
$e(\phi_{j,\pm}:i\lambda:x)$ as lying in $W_i\subset X$ and denote
such a section  by $E^i(\phi_{j,\pm}:i\lambda:x)$. Then, for
$\phi=(\phi_0,\phi_c)\in L^2(X-W_i,E)\oplus L^2(W_i,E)$ with
$\phi_c\in C^{\infty}_0(W_i,E)$ and $\phi_c|_{\partial W_i} \in
L$,
\begin{align}\label{e:E0}
\Dd_0(i)\phi \ = \  &
\frac{1}{4\pi}\sum_{j=1}^{d(\sigma_+)}\int^{\infty}_{-\infty}\lambda
E^i( \phi_{j,+}:i\lambda:x) \int_X E^i(
\phi_{j,+}:-i\lambda:y)\phi(y)
\ dy\ d\lambda \\
 - & \frac{1}{4\pi}\sum_{j=1}^{d(\sigma_+)}\int^{\infty}_{-\infty}
\lambda E^i( \phi_{j,-}:i\lambda:x) \int_X
E^i(\phi_{j,-}:-i\lambda:y)\phi(y)\ dy\ d\lambda. \notag
\end{align}
The following equality can be proved in the same way as in the
proof of Proposition \ref{p:maass},
\begin{equation}\label{e:eisen}
\sum^\kappa_{i=1}\int_{X_R}|E^i(\phi_{j,\pm}:i\lambda:x)|^2 dx =
2\kappa R + O(e^{-cR})
\end{equation}
for some positive constant $c$. (Let us remark that there is no
contribution over $\partial W_i$ by the choice of $\phi_{j,\pm}$
when we apply the Green formula as in \eqref{e:maass1}.) It
follows from (\ref{e:Eisen}), (\ref{e:eisen}) that
\begin{align}\label{e:trace-r}
& \int_{X} \vert E(\Phi_{j,\pm}:i\lambda:x)\vert^2- \sum^\kappa_{i=1}\vert E^i(\phi_{j,\pm}:i\lambda:x)\vert^2 dx\\
=&\lim_{R\to\infty} \
 \int_{X_R}  \vert E(\Phi_{j,\pm}:i\lambda:x)\vert^2- \sum^\kappa_{i=1}
 \vert E^i(\phi_{j,\pm}:i\lambda:x)\vert^2 dx \notag \\
=&-\tr(C_{\mp}(-i\lambda) C'_{\pm}(i\lambda)).\notag
\end{align}
Finally, the following proposition is the result of
(\ref{e:etrace}), (\ref{e:otrace}), (\ref{e:oid}),
(\ref{e:trace-e}), (\ref{e:trace-o}),  \eqref{e:E0},
(\ref{e:trace-r}) and Corollary \ref{c:ounip}.

\begin{prop} We have
\begin{align}\label{e:even0}
&\Tr(e^{-t\Dd^2}-\ski e^{-t\Dd^2_0(i)})\\
=&\sum_{\lambda_k} e^{-t\lambda^2_k}
-\frac{\dsi}{4\pi}\int^{\infty}_{-\infty} e^{-t\lambda^2}
\tr(C_{\Gamma}(\sigma_{+},-i\lambda)
C'_{\Gamma}(\sigma_{+},i\lambda))\
d\lambda \notag \\
=&\ I_{\G}(K^e_t)+H_{\G}(K^e_t)+U_{\G}(K^e_t), \notag
\end{align}
\begin{align}\label{e:odd}
&\Tr(\Dd e^{-t\Dd^2}-\ski \Dd_0(i) e^{-t\Dd^2_0(i)})\\
=&\sum_{\lambda_k} \lambda_k e^{-t\lambda^2_k}
-\frac{\dsi}{4\pi}\int^{\infty}_{-\infty} \lambda e^{-t\lambda^2}
\tr_s(C_{\Gamma}(\sigma_{+},-i\lambda)
C'_{\Gamma}(\sigma_{+},i\lambda))\
d\lambda \notag \\
=&\sum_{\lambda_k} \lambda_k e^{-t\lambda^2_k} -\frac{\dsi}{4\pi
}\int^{\infty}_{-\infty}\lambda e^{-t\lambda^2}
\tr(C_-(-i\lambda) C'_+(i\lambda))\ d\lambda \notag\\
&  \ \ \ \ \ \ \ \ \ \ \ \ \ \  +\frac{\dsi}{4\pi
}\int^{\infty}_{-\infty}\lambda e^{-t\lambda^2}\tr(C_+(-i\lambda)
C'_-(i\lambda))\
d\lambda \notag \\
=&\ H_{\G}(K^o_t).\notag
\end{align}

\end{prop}

\begin{rem}\em Although the definition of $\Dd_0(i)$ depends on the
choice of the Lagrangian subspace $L$, the relative traces on the
left sides of (\ref{e:even0}), (\ref{e:odd}) do not depend on this
choice. This is because the right side of \eqref{e:eisen} does not
depend on the choice of $L$.
\end{rem}

\section{Meromorphic continuations of the eta and zeta
functions}\label{s5}

 In this section, we prove Theorem \ref{t:pole}, which provides
us with the pole structures of the eta function and the zeta
function over $\bC$. We follow \cite{Mu3} and \cite{Mu4} and
define
\begin{equation}\label{def-eta}
\begin{split}
\eta_{\Dd,1}(z)&:=\frac{1}{\Gamma(\frac{z+1}{2})}\int^{1}_{0}
t^{\frac{z-1}{2}}\ \Tr(\Dd e^{-t\Dd^2}- \ski\Dd_0(i)
e^{-t\Dd^2_0(i)})\
dt,\\
\zeta_{\Dd^2,1}(z)&:=\frac{1}{\Gamma(z)}\int^{1}_{0} t^{z-1}\ [\
\Tr(e^{-t\Dd^2}-\ski e^{-t\Dd^2_0(i)})-h\ ]\ dt
\end{split}
\end{equation} for $\Re(z)\gg 0$ and
\begin{equation}\label{def-zeta}
\begin{split}
\eta_{\Dd,2}(z)&:=\frac{1}{\Gamma(\frac{z+1}{2})}\int^{\infty}_{1}\
t^{\frac{z-1}{2}} \ \Tr(\Dd e^{-t\Dd^2}-\ski \Dd_0(i)e^{-t\Dd^2_0(i)})\ dt,\\
\zeta_{\Dd^2,2}(z)&:=\frac{1}{\Gamma(z)}\int^{\infty}_{1} t^{z-1}\
[\ \Tr(e^{-t\Dd^2}-\ski e^{-t\Dd^2_0(i)})-h \ ]\ dt
\end{split}
\end{equation}
for $\Re(z)\ll 0$ where $h$ is the multiplicity of zero
eigenvalues of $\Dd^2$. To define the eta invariant and the
regularized determinant, we need to study the meromorphic
extensions of $\eta_{\Dd,i}(z)$ and $\zeta_{\Dd^2,i}(z)$ near
$z=0$ for $i=1,2$. The difficulty, which is not present for closed
manifolds, is the presence of continuous spectrum. Moreover, the
continuous spectrum of $\Dd$ is equal to the whole real line in
our case. Hence the meromorphic extensions of $\eta_{\Dd,2}(z)$
and $\zeta_{\Dd^2,2}(z)$ have nontrivial poles.

We start with $\eta_{\Dd,1}(z)$. It follows from (\ref{e:ozeta})
and (\ref{e:odd}) that
\begin{align}\label{e:odd-1}
&\Tr(\Dd e^{-t\Dd ^2}-\ski \Dd_0(i) e^{-t\Dd^2_0(i)})= H_{\G}(K^o_t)\\
=&\frac{2\pi i}{{(4\pi t)}^{\frac 32}}
\sum_{\gamma:\text{hyperbolic}}
l(C_{\gamma})^2j(\gamma)^{-1}D(\gamma)^{-1}\big(\overline{\chi_{\sigma_{+}}(m_{\gamma})}-
\overline{\chi_{\sigma_{-}}(m_{\gamma})}\big)\,
e^{-\frac{l(C_{\gamma})^2}{4t}}.\notag
\end{align}
The number
$c:=\text{min}_{\{\gamma:\text{hyperbolic}\}}l(C_{\gamma})$ is a
positive real number, hence as $t\to 0$,
\begin{equation}\label{e:asym-o-0}
\frac{2\pi i}{{(4\pi t)}^{\frac 32}}
\sum_{\gamma:\text{hyperbolic}}
l(C_{\gamma})^2j(\gamma)^{-1}D(\gamma)^{-1}\big(\overline{\chi_{\sigma_{+}}(m_{\gamma})}-
\overline{\chi_{\sigma_{-}}(m_{\gamma})}\big)\,
e^{-\frac{l(C_{\gamma})^2}{4t}} \ \sim\ a\, e^{-\frac{c^2}{4t}}
\end{equation}
for a constant $a$. Now (\ref{e:odd-1}) and (\ref{e:asym-o-0})
give
$$
\Tr(\Dd e^{-t\Dd^2}- \ski \Dd_0(i)e^{-t\Dd^2_0(i)})\ \sim \ a\,
e^{-\frac{c^2}{4t}} \quad\text{as}\quad t\to 0.
$$
This means that $\eta_{\Dd,1}(z)$ can be extended to the whole
complex plane without poles.

For $\zeta_{\Dd^2,1}(z)$, we use (\ref{e:even0}) to get
$$
\Tr(e^{-t\Dd^2}-\ski
e^{-t\Dd^2_0(i)})=I_{\G}(K^e_t)+H_{\G}(K^e_t)+U_{\G}(K^e_t).
$$
Recall that
\begin{equation}\label{e:ezeta1}
H_{\G}(K^e_t)=\frac{1}{\sqrt{4\pi t}}
\sum_{\gamma:\text{hyperbolic}}
l(C_{\gamma})j(\gamma)^{-1}D(\gamma)^{-1}\big(\overline{\chi_{\sigma_{+}}(m_{\gamma})}+
\overline{\chi_{\sigma_{-}}(m_{\gamma})}\big)\,
e^{-\frac{l(C_{\gamma})^2}{4t}}.
\end{equation}
As $H_{\G}(K^o_t)$, we have
\begin{equation}\label{e:hypc}
H_{\G}(K^e_t)\ \sim \ a'\, e^{-\frac{c^2}{4t}} \quad \text{as} \ \
t\to 0
\end{equation}
for a constant $a'$. We can see that $H_{\G}(K^e_t)$ does not give
any poles in the meromorphic continuation of $\zeta_{\Dd^2,1}(z)$.
Next we consider $I_{\G}(K^e_t)$ and $U_{\G}(K^e_t)$. An
elementary computation shows
\begin{equation}\label{e:idc}
I_{\G}(K^e_t)=\sum^{n}_{k=0} a_k t^{-k-\frac 12}
\end{equation}
for some constants $a_k$. By Corollary \ref{c:ounip}, we have
$$
U_{\G}(K^e_t)=\frac{2}{2\pi}\int^{\infty}_{-\infty}
e^{-t\lambda^2} (P_U(\lambda)+Q(\lambda))\ d\lambda
$$
where $P_U(\lambda)$ is an even polynomial of degree $(2n-4)$ and
$Q(\lambda)$ is given by
$$
{-\kappa \frac{d(\sigma_{+})}{2}}\Bigr(\psi(i\lambda+ \frac 12)
+\psi(-i\lambda +\frac 12) \Bigr).
$$
An elementary computation leads to
\begin{equation}\label{e:unipc1}
\int^{\infty}_{-\infty} e^{-t\lambda^2} P_U(\lambda)\
d\lambda=\sum^{n-2}_{k=0} b_k t^{-k-\frac 12}
\end{equation}
for some constants $b_k$. Using the relations
$$
\psi(x+1)=\frac{1}{x}+\psi(x), \qquad \psi(x)+\psi(x+\frac
12)=2(\psi(2x)-\log 2),
$$
we have
\begin{align*}
Q(\lambda)={- \frac{\kappa d(\sigma_{+})}{2}}\Big( \, 2\big(
\psi(2i\lambda+1)+\psi(-2i\lambda+1)-2\log 2 \big) - \big(
\psi(i\lambda+1)+\psi(-i\lambda+1) \big)\, \Big).
\end{align*}
To deal with the digamma function $\psi$, we use the following
asymptotic expansion
$$
\psi(z+1)\ \sim \ \log z+\frac{1}{2z}
-\sum^{\infty}_{k=1}\frac{B_{2k}}{(2k)z^{2k}} \qquad \text{as}
\quad z\to\infty
$$
where $B_{2k}$ are Bernoulli numbers. This implies the expansion
\begin{equation}\label{e:unipc3}
\int^{\infty}_{-\infty} e^{-t\lambda^2} Q(\lambda)\ d\lambda \
\sim \ \sum^{\infty}_{k=0} c_k t^{k-\frac 12} + d_0 t^{-\frac
12}\log t\quad \quad \text{as}\quad t\to 0
\end{equation}
for some constants $c_k$ and $d_0$.  By (\ref{e:hypc}),
(\ref{e:idc}), (\ref{e:unipc1}), (\ref{e:unipc3}), we have
$$
I_{\G}(K^e_t)+H_{\G}(K^e_t)+U_{\G}(K^e_t)\ \sim \
\sum^{\infty}_{k=-n} \beta_k t^{k-\frac 12} +\beta'_0 t^{-\frac
12}\log t\quad\quad \text{as} \quad t\to 0
$$
for constants $\beta_k, \beta'_0$. Therefore $\zeta_{\Dd^2,1}(z)$
is well defined for $\Re(z) > n+\frac 12$ and we can extend
$\zeta_{\Dd^2,1}(z)$ to a meromorphic function on $\bC$, with
poles determined by
\begin{equation}\label{e:zeta1}
\Gamma(z)\zeta_{\Dd^2,1}(z)=\sum^{\infty}_{k=-n}\frac{\beta_k}{z+k-\frac
12} +\frac{\beta'_0}{(z-\frac 12)^2} -\frac{h}{z}+H_1(z)
\end{equation}
where $\beta_k,\beta'_0$ are constants and $H_1(z)$ is a
holomorphic function.

To deal with the meromorphic extensions of $\eta_{\Dd,2}(z)$ and
$\zeta_{\Dd^2,2}(z)$, we consider the right sides of the following
equalities,
\begin{align*}
&\qquad \Tr(\Dd e^{-t\Dd^2}-\ski \Dd_0(i) e^{-t\Dd^2_0(i)})=
\sum_{\lambda_k} \lambda_k
e^{-t\lambda_k}\\
&\qquad\qquad\qquad -\frac{\dsi}{4\pi}\int^{\infty}_{-\infty}
\lambda e^{-t\lambda^2} \tr_s(C_{\Gamma}(\sigma_{+},-i\lambda)
C'_{\Gamma}(\sigma_{+},i\lambda))\ d\lambda,
\end{align*}
\begin{align*}
&\Tr(e^{-t\Dd^2}- \ski e^{-t\Dd^2_0(i)})-h= \sum_{\lambda_k\neq 0}
e^{-t\lambda^2_k}\\
&\qquad -\frac{\dsi}{4\pi}\int^{\infty}_{-\infty} e^{-t\lambda^2}
\tr(C_{\Gamma}(\sigma_{+},-i\lambda)
C'_{\Gamma}(\sigma_{+},i\lambda))\ d\lambda.
\end{align*}
The discrete eigenvalues give
$$
\sum_{\lambda_k} \lambda_k e^{-t\lambda^2_k}\ \sim\ e^{-ct}, \quad
\sum_{\lambda_k\neq 0}e^{-t\lambda^2_k} \  \sim \
e^{-ct}\qquad\text{as}\ \ t\to\infty
$$
for a positive constant $c$. The operator
$C_{\Gamma}(\sigma_{+},i\lambda)$ is analytic along the imaginary
axis and
\begin{align*}
\tr_s(C_{\Gamma}(\sigma_{+},-i\lambda)
C'_{\Gamma}(\sigma_{+},i\lambda))&=
-\tr_s(C_{\Gamma}(\sigma_{+},i\lambda) C'_{\Gamma}(\sigma_{+},-i\lambda)),\\
\tr(C_{\Gamma}(\sigma_{+},-i\lambda)
C'_{\Gamma}(\sigma_{+},i\lambda))&=
\tr(C_{\Gamma}(\sigma_{+},i\lambda)
C'_{\Gamma}(\sigma_{+},-i\lambda)),
\end{align*}
hence we have the following analytic expansion at $\lambda=0$:
$$
\tr_s(C_{\Gamma}(\sigma_{+},-i\lambda)
C'_{\Gamma}(\sigma_{+},i\lambda)) =\sum^{\infty}_{k=0} f_{2k+1}
\lambda^{2k+1},
$$
$$
\tr(C_{\Gamma}(\sigma_{+},-i\lambda)
C'_{\Gamma}(\sigma_{+},i\lambda)) =\sum^{\infty}_{k=0}
g_{2k}\lambda^{2k}
$$
for some constants $f_{2k+1}, g_{2k}$. Therefore we have
\begin{equation}\label{e:ol}
\int^{1}_{-1} \lambda e^{-t\lambda^2}
\tr_s(C_{\Gamma}(\sigma_{+},-i\lambda)
C'_{\Gamma}(\sigma_{+},i\lambda))\ d\lambda  \quad \sim \quad
\sum^{\infty}_{k=0} \gamma_k t^{-(k+\frac 32)} \quad \text{as} \
t\to\infty,
\end{equation}
\begin{equation}\label{e:el}
\int^{1}_{-1}  e^{-t\lambda^2}
\tr(C_{\Gamma}(\sigma_{+},-i\lambda)
C'_{\Gamma}(\sigma_{+},i\lambda))\ d\lambda \quad \sim \quad
\sum^{\infty}_{k=0} \gamma'_k t^{-(k+\frac{1}{2})} \quad \text{as}
\ t\to\infty.
\end{equation}
The corresponding integrals over $(-\infty,1]_{\lambda}\cup
[1,\infty)_\lambda$ converge to $0$ exponentially as $t\to\infty$.
Hence, the expansion (\ref{e:ol}) shows that $\eta_{\Dd,2}(z)$ is
well defined for $\Re(z)< 1 $ and extends to the whole complex
plane with the following pole structure:
\begin{equation}\label{e:eta2}
\Gamma((z+1)/2)\eta_{\Dd,2}(z)=\sum^{\infty}_{k=0}\frac{-2\gamma_k}{z-2k-2}+K_2(z)
\end{equation}
for constants $\gamma_k$ and a holomorphic function $K_2(z)$. In
the same way, the expansion (\ref{e:el}) implies that
$\zeta_{\Dd^2,2}(z)$ is well defined for $\Re(z)< \frac 12$ and
extends to the whole complex plane with poles determined by the
equality
\begin{equation}\label{e:zeta2}
\Gamma(z)\zeta_{\Dd^2,2}(z)=\sum^{\infty}_{k=0}\frac{-\gamma'_k}{z-k-\frac
12}+H_2(z)
\end{equation}
for the constants $\gamma'_k$ and a holomorphic function $H_2(z)$.

We define the eta and zeta functions by
\begin{equation}\label{eta-zeta}
\eta_{\Dd}(z):=\eta_{\Dd,1}(z)+\eta_{\Dd,2}(z), \quad
\zeta_{\Dd^2}(z):=\zeta_{\Dd^2,1}(z)+\zeta_{\Dd^2,2}(z).
\end{equation}
Here now the right sides of these equalities are meromorphic
functions over $\bC$ with the poles described in the above. The
equalities (\ref{e:zeta1}), (\ref{e:eta2}), (\ref{e:zeta2}) give
the following theorem

\begin{thm} The poles of the eta function $\eta_{\Dd}(z)$ and
the zeta function $\zeta_{\Dd^2}(z)$ are determined by the
equations
\begin{align*}
&\Gamma((z+1)/2)\eta_{\Dd}(z)=\sum^{\infty}_{k=0}\frac{-2\gamma_k}{z-2k-2}+K(z),\\
&\Gamma(z)\zeta_{\Dd^2}(z)=\sum^{\infty}_{k=-n}\frac{\beta_k}{z+k-\frac
12} +\frac{\beta'_0}{(z-\frac 12)^2}-\frac{h}{z}
+\sum^{\infty}_{k=0}\frac{-\gamma'_k}{z-k-\frac 12}+H(z)
\end{align*}
where $K(z)$ and $H(z)$ are holomorphic. In particular,
$\eta_{\Dd}(z)$ and $\zeta_{\Dd^2}(z)$ are regular at $z=0$.
\end{thm}

\section{Eta invariants, Zeta functions of odd type and functional equations}
\label{s6}

In this section, we study the eta invariant and its relation with
the Selberg zeta function of odd type. We use the Selberg trace
formula to prove a generalization of Millson's theorem in \cite
{M} for hyperbolic manifolds with cusps. Since the unipotent term
in our situation vanishes (see Corollary \ref{c:ounip}), we obtain
the same formula as in the case of Millson in \cite{M}. We also
derive the functional equation for the eta invariant and the
Selberg zeta function of odd type.

By the analysis in Section \ref{s5}, the eta function
$\eta_{\Dd}(z)$ is regular at $z=0$ and we can put $z=0$ in the
equality
$$
\eta_{\Dd}(z) \ = \ \frac{1}{\Gamma(\frac{z+1}{2})}\int^{\infty}_0
t^{\frac{z- 1}{2}}  \Tr(\Dd e^{-t\Dd^2}-\ski \Dd_0(i)
e^{-t\Dd^2_0(i)})\ dt.
$$
We define the eta invariant of $\Dd$ by
$$
\eta(\Dd):=\ \eta_{\Dd}(0) \ = \
\frac{1}{\sqrt{\pi}}\int^{\infty}_0 t^{-\frac 12}  \Tr(\Dd
e^{-t\Dd^2}- \ski \Dd_0(i) e^{-t\Dd^2_0(i)})\ dt.
$$
Let us recall that \begin{multline*} \Tr(\Dd e^{-t\Dd^2}-\ski
\Dd_0(i) e^{-t\Dd^2_0(i)}) \\=\frac{2\pi i}{{(4\pi t)}^{\frac 32}}
\sum_{\gamma:\text{hyperbolic}}
l(C_{\gamma})^2j(\gamma)^{-1}D(\gamma)^{-1}\big(\overline{\chi_{\sigma_{+}}(m_{\gamma})}-
\overline{\chi_{\sigma_{-}}(m_{\gamma})}\big)\,
e^{-\frac{l(C_{\gamma})^2}{4t}}.
\end{multline*}
Using the elementary equality
$\int^{\infty}_0e^{-s^2t}\frac{e^{-r^2/4t}}{(4\pi
t)^{\frac{3}{2}}}dt=\frac{e^{-sr}}{4\pi r}$, we have
\begin{align*}
&\int^{\infty}_{0}e^{-s^2t}\Tr(\Dd e^{-t\Dd^2}-\ski \Dd_0(i)
e^{-t\Dd^2_0(i)})\ dt\\
= &\frac{i}{2}\sum_{\gamma:\text{hyperbolic}}
l(C_{\gamma})j(\gamma)^{-1}D(\gamma)^{-1}\big(\overline{\chi_{\sigma_{+}}(m_{\gamma})}-
\overline{\chi_{\sigma_{-}}(m_{\gamma})}\big)\, e^{-sl(C_{\gamma})}\\
= &\frac{i}{2}\sum_{\gamma:\text{hyperbolic}}
l(C_{\gamma})j(\gamma)^{-1}\vert
\text{det}(\text{Ad}(a_{\gamma}m_{\gamma})^{-1} -
I\vert_{\frak{n}})\vert^{-1}\big(\overline{\chi_{\sigma_{+}}(m_{\gamma})}-
\overline{\chi_{\sigma_{-}}(m_{\gamma})}\big)\,
e^{-(s+n)l(C_{\gamma})}.
\end{align*}
We define the Selberg zeta function of odd type by
\begin{multline}\label{def-odd}
Z^o_H(s)\\
:=\exp\big(-\sum_{\gamma:\text{hyperbolic}}j(\gamma)^{-1}\vert
\text{det}(\text{Ad}(a_{\gamma}m_{\gamma})^{-1} -
I\vert_{\frak{n}})\vert^{-1}(\overline{\chi_{\sigma_{+}}(m_{\gamma})}-
\overline{\chi_{\sigma_{-}}(m_{\gamma})})e^{-sl(C_{\gamma})}\big)
\end{multline}
for $\Re(s)\gg 0$. In Proposition \ref{p:reg} we will show that
$Z^o_H(s)$ has a meromorphic extension over $\bC$ and $Z^o_H(s)$ is
regular at $s=n$. Now we have
\begin{equation*}\label{e:ortr}
\int^{\infty}_{0}e^{-s^2t}\Tr(\Dd e^{-t\Dd^2}-\ski\Dd_0(i)
e^{-t\Dd^2_0(i)})\ dt =\frac{i}{2}\frac{d}{ds}\log Z^o_H(s+n).
\end{equation*}
Following the argument on p.\ 27 of \cite{M}, we use the equality
$$
t^{\frac{z-1}2} \ = \ \frac{2}{\Gamma(\frac{1-z}{2})}
\int^{\infty}_0 s^{-z} e^{-s^2t} \ ds
$$
for $\Re(z)<1$ to get
\begin{align*}
\eta_{\Dd}(z) \ =& \
\frac{2}{\Gamma(\frac{z+1}{2})\Gamma(\frac{1-z}{2})}
\int^{\infty}_0 \int^{\infty}_0 s^{-z} e^{-s^2t} \ ds\ \Tr(\Dd
e^{-t\Dd^2}- \ski \Dd_0(i)
e^{-t\Dd^2_0(i)})\ dt \\
=& \ \frac{2}{\Gamma(\frac{z+1}{2})\Gamma(\frac{1-z}{2})}
\int^{\infty}_0 s^{-z} \int^{\infty}_0  e^{-s^2t}   \Tr(\Dd
e^{-t\Dd^2}-\ski \Dd_0(i)
e^{-t\Dd^2_0(i)}) \ dt \ ds \\
=& \ \frac{i}{\Gamma(\frac{z+1}{2})\Gamma(\frac{1-z}{2})}
\int^{\infty}_0 s^{-z}  \frac{d}{ds} \log Z^o_H(s+n)  \ ds .
\end{align*}
If we evaluate the above equality at $z=0$, we get the following
theorem.

\begin{thm} \label{t:eta}
For the Dirac operator $\Dd$ over a $(2n+1)$-dimensional
hyperbolic manifold with cusps, the eta invariant $\eta(\Dd)$ and
the Selberg zeta function of odd type $Z^o_H(s)$ satisfy
\begin{equation}\label{e:eta}
\eta(\Dd)=\frac{1}{\pi i}\log Z^o_H(n).
\end{equation}

\end{thm}

\vspace{0.3cm}

Now let us show that $Z^o_H(s)$ has the meromorphic extension over
$\bC$. We select a smooth odd function $g(u)$ such that $|g(u)|=1$
if $|u| > c$, $g(u)=0$ near $0$ and $\int^{\infty}_0 g'(u) du=1$
where $c=\text{min}_{\{\gamma:\text{hyperbolic}\}}l(C_{\gamma})$.
We define
$$
\left\{ \aligned &H_s(\sigma_+,\lambda):=\int^{\infty}_{-\infty}
g(u)e^{-s|u|}e^{i\lambda u}du, \\
                 &H_s(\sigma_-,\lambda):=-H_s(\sigma_+,\lambda)
        \endaligned \right.
$$
for a complex parameter $s$. By the Paley-Wiener theorem, ( see
theorem 2.2 in \cite{M} ), there exists $f_s$ over $G$ with
$\hat{f}_s(\sigma_{\pm},i\lambda)=H_s(\sigma_{\pm}, \lambda)$,
$\hat{f}_s(\sigma,i\lambda)=0$ if $\sigma\neq \sigma_{\pm}$.
Applying the Selberg trace formula to the one parameter family of
functions $f_{s}$ on $G$ for $\Re(s)\gg 0$, we get
\begin{align}\label{e:odd-sca}
&\sum_{\lambda_j\in \sigma^+_p}H_s(\sigma_+,\lambda_j)
+\sum_{\lambda_j\in \sigma^-_p}H_s(\sigma_-,\lambda_j) \\
&-\frac{\dsi}{4\pi}\int^{\infty}_{-\infty}H_s(\sigma_+,\lambda)
\tr_s(C_{\Gamma}(\sigma_{+},-i\lambda)
C'_{\Gamma}(\sigma_{+},i\lambda))\ d\lambda\notag \\
=&\sum_{\gamma:\text{hyperbolic}}
l(C_{\gamma})j(\gamma)^{-1}D(\gamma)^{-1}\big(\overline{\chi_{\sigma_{+}}(m_{\gamma})}-
\overline{\chi_{\sigma_{-}}(m_{\gamma})}\big)\, e^{-sl(C_{\gamma})} \notag \\
=&\ \frac{d}{ds}\log Z^o_H(s+n). \notag
\end{align}
Here we used the fact that the identity, unipotent orbital
integrals vanish by the definition of $H_s(\sigma_{\pm},\lambda)$
and results in the previous sections. We shall use this equality
to get the meromorphic extension of
$$
\mathcal{Z}(s):=\frac{d}{ds}\log Z^o_H(s+n)
$$
over $\bC$ and to investigate its poles.

\textbf{Discrete eigenvalue term:} \ Integration by part gives
\begin{equation*}
H_s(\sigma_{\pm},\lambda)= \pm \frac{1}{s-i\lambda}\int^{\infty}_0
g'(u)e^{-(s-i\lambda)u}du  \mp
\frac{1}{s+i\lambda}\int^{\infty}_{0} g'(u)e^{-(s+i\lambda)u}du.
\end{equation*}
This equality provides the meromorphic extension of
$H_s(\sigma_{\pm},\lambda_j)$ over $\bC$ and we see that
$H_s(\sigma_{\pm},\lambda_j)$ has simple poles at $i\lambda_j$,
$-i\lambda_j$ with residues $\pm m_j$, $\mp m_j$ for $\lambda_j\in
\sigma^{\pm}_p$ where $m_j$ is the multiplicity of $\lambda_j$.

\textbf{Scattering term:} \ We consider the scattering term
\begin{equation}\label{e:last}
-\frac{\dsi}{4\pi}\int^{\infty}_{-\infty}H_s(\sigma_+,\lambda)
\tr_s(C_{\Gamma}(\sigma_{+},-i\lambda)
C'_{\Gamma}(\sigma_{+},i\lambda))\ d\lambda.
\end{equation}
First we observe
\begin{align*}
\Phi(z):=& \ \tr_s(C_{\Gamma}(\sigma_{+},-z)
\frac{d}{dz}C_{\Gamma}(\sigma_{+},z)) \ = \ \frac{d}{dz} \biggl(\
\log\det C_+(z) - \log\det C_-(z) \ \biggr).
\end{align*}
Using the equalities
$C_{\G}(\sigma_+,z)C_{\G}(\sigma_+,-{z})=\Id$,
$C_{\Gamma}(\sigma_+,z)^*=C_{\Gamma}(\sigma_+,\bar{z})$ and the
first displayed formula on p. 518 and (6.8) of \cite{Mu89}, we
have
$$
\det C_+(z) \ = \det C_+(0)\ p_+^z\prod_{\Re(q_+)<0}
\frac{z+\bar{q}_+}{z-q_+},
$$
$$
\det C_-(z) \ = \det C_-(0)\ p_-^{z}\prod_{\Re(q_-)<0}
\frac{z+\bar{q}_-}{z-q_-}
$$
for some constants $p_+,p_-$. Here the infinite products are taken
over the poles $\{q_{\pm}\}$ of $\det C_\pm(z)$ respectively. Note
that $\det C_\pm(z)$ is holomorphic over the half plane with
$\Re(z)\geq 0$. Hence $\Phi(z)$ has the following form over $\bC$:
\begin{align}\label{e:infs}
\Phi(z)  =  -\sum_{\Re(q_+)<0}\frac{2\,
\Re(q_+)}{(z-q_+)(z+\bar{q}_+)}  + \sum_{\Re(q_-)<0} \frac{2\,
\Re(q_-)}{(z-q_-)(z+\bar{q}_-)}  + \log \ p_+ - \log \ p_-.
\end{align}
Now we consider the contour integral
$$
\mathcal{L}_R:=\frac{1}{4\pi i}\int_{L_R} H_s(\sigma_+, z) \
\tr_s(C_{\Gamma}(\sigma_{+},-iz) C'_{\Gamma}(\sigma_{+},iz)) \, dz
$$
where $L_R=[\ -R,\ R \ ]\ \cup \{\ Re^{i\theta} \ | \ 0\le
\theta\le \pi\ \}$. As in proposition 3.10 of \cite{GW}, we can
apply the Cauchy integral formula and obtain
\begin{multline}\label{e:mero-scattering}
\lim_{R\to\infty}\mathcal{L}_{R} \\ = i \Big( \sum_{\Re(q_+)<0}
\frac{1}{s-q_+}\int^{\infty}_0 g'(u)e^{-(s-q_+)u}du  -
\sum_{\Re(q_-)<0} \frac{1}{s-q_-}\int^{\infty}_0
g'(u)e^{-(s-q_-)u}du \Big)
\end{multline}
for $\Re(s) \gg 0$.  Now the right side of
\eqref{e:mero-scattering} gives us the meromorphic extension over
$\bC$. On the other hand, we can also show that
$$
\lim_{R\to\infty} \mathcal{L}_R =\frac{1}{4\pi i}
\int^{\infty}_{-\infty} H_s(\sigma_+,\lambda)
\tr_s(C_{\Gamma}(\sigma_{+},-i\lambda)
C'_{\Gamma}(\sigma_{+},i\lambda))\, d\lambda.
$$
Here the integral over the semicircle of radius $R$ vanishes as
$R\to\infty$ by the definition of $H_s(\sigma_+,z)$ and
\eqref{e:infs}. Therefore the meromorphic extension of the
scattering term has simple poles at $q_{\pm}$ for $\Re(q_{\pm})<0$
with residues $\pm d(\sigma_+)b_{\pm}$ where $b_{\pm}$ denotes the
order of the pole of $\det C_{\G}(\sigma_+,z)$ at $q_{\pm}$.

 Combining the contributions of the
$H_s(\sigma_{+},\lambda_j)$'s and the scattering term, we see that
$\mathcal{Z}(s)$ has a meromorphic extension over $\bC$ and its
simple poles are located at $\pm i\lambda_j$, $\mp i\lambda_j$ for
$\lambda_j\in \sigma^{\pm}_p$, and at $q_{\pm}$ for
$\Re(q_{\pm})<0$ with residues $m_j$, $-m_j$ and
$\pm\dsi{b_{\pm}}$ respectively. In particular, we can see that
$\mathcal{Z}(s)$ is regular at $s=0$. This implies the following
proposition.

\begin{prop}\label{p:reg} The Selberg zeta function of odd type
$Z^o_H(s+n)$ has a meromorphic extension over $\bC$, and is
regular at $s=0$.
\end{prop}

\begin{rem}\em
The zeros of $Z^o_H(s+n)$ are located at $\pm i\lambda_j$ for
$\lambda_j\in \sigma^{\pm}_p$, at $q_+$ for $\Re(q_+)<0$ and their
orders are $m_j$, $d(\sigma_+)b_+$. The  poles of $Z^o_H(s+n)$ are
located at $\mp i\lambda_j$ for $\lambda_j\in\sigma^{\pm}_p$, at
$q_-$ for $\Re(q_-)<0$ and their orders are $m_j$,
$d(\sigma_+)b_-$.
\end{rem}

Let us study the functional equation of $\eta(\Dd)$ and
$Z^o_H(s)$. We set
$$
R(s):=\mathcal{Z}(s) -\mathcal{Z}(-s)+ \dsi\Phi(s).
$$
Then $R(s)$ is an odd entire function of $s$. Let $h(s)$ be an odd
function which decreases sufficiently rapidly as
$\text{Im}(s)\to\infty$ in the strip $ \{ \ s\in\mathbb{C} \ | \
|\Re(s)| < n+\epsilon, \ \epsilon>0\ \}$ and consider the contour
integral
$$
\mathcal{L}_T:=\frac{1}{2\pi i}\int_{L_T} h(s)\mathcal {Z}(s)\, ds
$$
where $L_T$ is the rectangle with the corners $ a+ iT, a-iT,
-a+iT, -a-iT$ with $n<a <n+\epsilon$. Then we have
$$
\lim_{T\to\infty}\mathcal{L}_T =\frac{1}{2\pi
i}\int^{a+i\infty}_{a-i\infty}h(s)(\mathcal{Z}(s)
-\mathcal{Z}(-s))\, ds +2\frac{1}{2\pi
i}\int^{-a-i\infty}_{-a+i\infty}h(s) \mathcal{Z}(s)\, ds.
$$
We apply the Cauchy integral theorem to get the equality
\begin{equation}\label{e:LT}
\lim_{T\to\infty}\mathcal{L}_T= 2 \sum_{\lambda_j}m_jh(i\lambda_j)
+\sum_{-a< \Re(q_k)< 0} \dsi b_k h(q_k)
\end{equation}
where we use the notations $q_k,b_k$ instead of $q_{\pm}$, $\pm
b_{\pm}$. Because the simple poles of $\mathcal{Z}(s)$ are located
at the $q_k$'s with residues $\dsi b_k$, between $\Re(s)=-a$ and
$\Re(s)=0$, we have
\begin{equation}\label{e:right}
\frac{1}{2\pi i}\int^{-a-i\infty}_{-a+i\infty}h(s)
\mathcal{Z}(s)\, ds =\frac{1}{2\pi
i}\int^{-i\infty}_{i\infty}h(\lambda) \mathcal{Z}(\lambda) \
d\lambda +\sum_{-a< \Re(q_k)< 0} \dsi b_k h(q_k).
\end{equation}
The simple poles at $-q_k$ of $\mathcal{Z}(-s)$ between $\Re(s)=a$
and $\Re(s)=0$ give rise to the equality
\begin{align}\label{e:R(s)}
&\frac{1}{2\pi i}\int^{a+i\infty}_{a-i\infty}h(s)(\mathcal{Z}(s)
-\mathcal{Z}(-s))\, ds \\
=&\frac{1}{2\pi
i}\int^{i\infty}_{-i\infty}h(\lambda)(\mathcal{Z}(\lambda)
 -\mathcal{Z}(-\lambda))\, d\lambda
+\sum_{0< \Re(-q_k)< a} \dsi b_k h(-q_k).\notag
\end{align}
By (\ref{e:LT}), (\ref{e:right}) and (\ref{e:R(s)}), we have
\begin{align}\label{e:second}
2\sum_{\lambda_j}m_jh(i\lambda_j) =&\  2\frac{1}{2\pi
i}\int^{-i\infty}_{i\infty}h(\lambda)
(\mathcal{Z}(\lambda) + \frac{\dsi}{2}\Phi(\lambda))\ d\lambda\\
& +\frac{1}{2\pi
i}\int^{i\infty}_{-i\infty}h(\lambda)(\mathcal{Z}(\lambda)
-\mathcal{Z}(-\lambda) + \dsi\Phi(\lambda))\ d\lambda.\notag
\end{align}
If we change variables $\lambda\to i\lambda$, we see that the first
term of the right side in (\ref{e:second}) is equal to two times of
\begin{align}\label{e:selberg}
&\ \frac{1}{2\pi}\int^{\infty}_{-\infty}h(i\lambda)
\sum_{\gamma:\text{hyperbolic}}
l(C_{\gamma})j(\gamma)^{-1}D(\gamma)^{-1}\big(\overline{\chi_{\sigma_{+}}(m_{\gamma})}-
\overline{\chi_{\sigma_{-}}(m_{\gamma})}\big)\, e^{-i\lambda
l(C_{\gamma})}\ d\lambda\\
&\ +\frac{\dsi}{4\pi }\int^{\infty}_{-\infty}h(i\lambda)
\Phi(i\lambda)\ d\lambda\notag\\
=&\  \sum_{\gamma:\text{hyperbolic}}
l(C_{\gamma})j(\gamma)^{-1}D(\gamma)^{-1}\big(\overline{\chi_{\sigma_{+}}(m_{\gamma})}-
\overline{\chi_{\sigma_{-}}(m_{\gamma})}\big)
\frac{1}{2\pi}\int^{\infty}_{-\infty} h(i\lambda) e^{-i\lambda l(C_{\gamma})} d\lambda \notag \\
&\  +\frac{\dsi}{4\pi }\int^{\infty}_{-\infty}h(i\lambda)
\Phi(i\lambda)\ d\lambda \ .\notag
\end{align}
By (\ref{e:second}), (\ref{e:selberg}) and the Selberg trace
formula applied to $h(i\lambda)$, we have
$$
\frac{1}{2\pi i}\int^{i\infty}_{-i\infty}h(\lambda)R(\lambda)\
d\lambda =0
$$
for any odd holomorphic $h(s)$ satisfying suitable growth
conditions. Since $R(s)$ is an entire function, this implies that
$R(s)=0$ over $\bC$. Hence we have the equality
$$
\frac{d}{ds}\log(Z^o_H(n+s)Z^o_H(n-s))\ =\ -\dsi \Phi(s)
$$
for any $s\in \bC$. Finally, recalling the equalities
$$
\Phi(s) \ =\ \frac{d}{ds} \biggl(\ \log\det C_+(s) - \log\det
C_-(s) \ \biggr),
$$
$$
\eta(\Dd)=\frac{1}{\pi i}\log Z^o_H(n),
$$
we get the following theorem.

\begin{thm}
We have
$$
Z^o_H(n+s)Z^o_H(n-s)=\exp(2\pi i\eta(\Dd))\ \biggl(\frac{\det
C_+(s) C_-(0)}{\det C_-(s) C_+(0)}\biggr)^{-\dsi}  \quad
\text{for} \quad s\in\bC .
$$
\end{thm}

\begin{rem} \em In the above equality, $Z^o_H(n+s)Z^o_H(n-s)$ has
zeros at $q_+,-q_+$ of orders $d(\sigma_+)b_+$ for $\Re(q_+)<0$
and poles at $q_-,-q_-$ of orders $d(\sigma_+)b_-$ for
$\Re(q_-)<0$. These zeros and poles coincide with the zeros and
poles of $(\frac{\det C_+(s) C_-(0)}{\det C_-(s)
C_+(0)})^{-\dsi}$.
\end{rem}

\section{Cusp contributions for regularized determinants and functional equations}
\label{s7}

In this section, we compute the unipotent factor in the relation
between the regularized determinant and the Selberg zeta function
of even type. We also derive the functional equation for the
regularized determinant and the Selberg zeta function of even type
where the unipotent factor plays a nontrivial role. This equation
is the even type counterpart of the functional equation for the
eta invariant and the Selberg zeta function of odd type proved in
Section \ref{s6}.

We define the zeta function with the factor $e^{-ts^2}$ for a
positive real number $s$ by
$$
\zeta_{\Dd^2}(z,s):=\frac{1}{\Gamma(z)}\int^{\infty}_{0} t^{z-1} \
\Tr(e^{-t\Dd^2}-\ski e^{-t\Dd^2_0(i)})\ e^{-ts^2} dt.
$$
Note that $\zeta_{\Dd^2}(z,0)=\zeta_{\Dd^2}(z)$ if the kernel of
$\Dd^2$ is trivial. We proved that $\zeta_{\Dd^2}(z)$ is regular
at $z=0$ in Section \ref{s5}. Due to the factor $s^2$, the
continuous spectrum of $\Dd^2+s^2$ has minimum $s^2$ and the large
time part of the zeta function  $\zeta_{\Dd^2}(z,s)$ does not
produce any poles. Therefore $\zeta_{\Dd^2}(z,s)$ is regular at
$z=0$. We define the regularized determinant of $\Dd^2+s^2$ by
$$
\Det(\Dd^2, s):=\exp(-\zeta'_{\Dd^2}(0,s))
$$
for a positive real number $s$. Now we observe that
$\zeta_{\Dd^2}(0,s)=0$ by the results of Section \ref{s5}. Then we
have
$$
\frac{d}{dz}\biggm|_{z=0}\zeta_{\Dd^2}(z,s)\ =\ \int^{\infty}_{0}
t^{z-1} \ \Tr(e^{-t\Dd^2}- \ski e^{-t\Dd^2_0(i)})\ e^{-ts^2} dt
\biggm|_{z=0} .
$$
Recall that
\begin{align}\label{e:even}
&\Tr(e^{-t\Dd^2}- \ski e^{-t\Dd^2_0(i)}) \\
=&\sum_{\lambda_k} e^{-t\lambda^2_k} -\frac{\dsi}{4\pi}
\int^{\infty}_{-\infty} e^{-t\lambda^2}
\tr(C_{\Gamma}(\sigma_{+},-i\lambda)
C'_{\Gamma}(\sigma_{+},i\lambda))\
d\lambda \notag \\
=&I_{\G}(K^e_t)+H_{\G}(K^e_t)+U_{\G}(K^e_t), \notag
\end{align}
and put
\begin{align*}
I(z,s)\ :=&\ \int^{\infty}_{0} t^{z-1} I_{\G}(K^e_t)\,  e^{-ts^2} dt,\\
H(z,s)\ :=&\ \int^{\infty}_{0} t^{z-1} H_{\G}(K^e_t)\, e^{-ts^2} dt,\\
U(z,s)\ :=&\ \int^{\infty}_{0} t^{z-1} U_{\G}(K^e_t)\, e^{-ts^2}
dt.
\end{align*}
Then we have
\begin{equation}\label{e:00}
\frac{1}{2s}\frac{d}{ds} \log \Det(\Dd^2,s) =
I(1,s)+H(1,s)+U(1,s).
\end{equation}
Now we want to find $Z_I(s), Z^e_H(s), Z_U(s)$ such that
$$
I(1,s)=\frac{1}{2s}\frac{d}{ds}\log Z_I(s), \ \
H(1,s)=\frac{1}{2s}\frac{d}{ds}\log Z^e_H(s), \ \
U(1,s)=\frac{1}{2s}\frac{d}{ds}\log Z_U(s).
$$

First,we recall that the Selberg zeta function of even type,
\begin{equation}\label{def-even}
Z^e_H(s):=\exp(-\sum_{\sigma_{\pm}}\sum_{\gamma:\text{hyperbolic}}
j(\gamma)^{-1}\vert
\text{det}(\text{Ad}(a_{\gamma}m_{\gamma})^{-1} -
I\vert_{\frak{n}})\vert^{-1}\overline{\chi_{\sigma}(m_{\gamma})}\,
e^{-sl(C_{\gamma})})
\end{equation}
is defined for $\Re(s)\gg 0$ and we will show that this has a
meromorphic extension over $\bC$. We use (\ref{e:ezeta}) and
elementary equality $\int^{\infty}_0 \frac{1}{\sqrt{4\pi
t}}\exp(-(\frac{r^2}{4t}+ts^2))dt= \frac{1}{2s} e^{-rs}$ to see
\begin{equation}\label{hypf}
H(1,s)\ =\ \frac{1}{2s}\frac{d}{ds}\log Z^e_H(s+n).
\end{equation}

It follows from Corollary \ref{c:ounip} that $U(z,s)$ is the
Mellin transform with factor $e^{-ts^2}$ of the following terms
\begin{equation} \label{U12}
\frac{2}{2\pi}\int^{\infty}_{-\infty} e^{-t\lambda^2}\Bigr (
P_U(\lambda)- {\kappa
\frac{d(\sigma_{\pm})}{2}}\big(\psi(i\lambda+ \frac 12)
+\psi(-i\lambda +\frac 12)\big) \Bigr )\ d\lambda
\end{equation}
where $P_U(\lambda)$ is an even polynomial of degree $(2n-4)$. We
deal with the term $\psi(\pm i\lambda+\frac 12)$ using the Cauchy
integral formula
\begin{align*}
&\frac{1}{2\pi}\int^{\infty}_0\int^{\infty}_{-\infty}
e^{-t(\lambda^2+s^2)}
\psi(\pm i\lambda+\frac 12)\ d\lambda \ dt\\
=&\frac{1}{2\pi}\int^{\infty}_{-\infty}
\frac{1}{\lambda^2+s^2}\psi(\pm i\lambda+\frac 12)\ d\lambda
=\frac{1}{2s}\psi(s+\frac 12).
\end{align*}
For the part $P_U(\lambda)$ in \eqref{U12},  we have
$$
\frac{1}{2s}\frac{d}{ds}\log \, \exp(\int^s_0 P_U(i\lambda)\
d\lambda) =\frac{1}{2\pi}\int^{\infty}_0
e^{-ts^2}\int^{\infty}_{-\infty} e^{-t\lambda^2}P_U(\lambda)\
d\lambda\ dt
$$
by lemma 3 in \cite{F}. We define
\[
Z_U(s):=\Gamma(s+\frac 12)^{-2\kappa d(\sigma_{+})} \exp(2\int^s_0
P_U(i\lambda)\ d\lambda). \] Then this satisfies
\begin{equation}\label{unipf}
U(1,s)=\frac{1}{2s}\frac{d}{ds}\log Z_U(s).
\end{equation}

For $I(z,s)$, we can treat this as for $P_U(\lambda)$ since
$p(\sigma_{+},\lambda)=p(\sigma_{-},\lambda)$ is an even
polynomial of $\lambda$. Hence, we  can see
\begin{equation}\label{e:idenf}
Z_I(s):=\exp(2\int^s_0 P_I(i\lambda)\ d\lambda)
\end{equation}
where $P_I(\lambda)={2\pi \mathrm{Vol}(\Gamma\bs {G})}\,
p(\sigma_{+},\lambda)$. By \eqref{e:00}, \eqref{hypf},
\eqref{unipf} and \eqref{e:idenf}, we have the following theorem

\begin{thm}\label{t:det} The following equality holds for any $s
\in \mathbb{C}$,
\begin{align}\label{e:det}
\Det(\Dd^2,s) \ =&\ C\, Z^e_H(s+n) \Gamma(s+\frac 12)^{-2\kappa
d(\sigma_{+})} \exp(2\int^s_0 P_I(i\lambda)+P_U(i\lambda)\
d\lambda)
\end{align}
where $C$ is a constant and $P_I(\lambda)={2\pi
\mathrm{Vol}(\Gamma\bs {G})}p(\sigma_{+},\lambda)$.
\end{thm}

\begin{proof}
A priori, the equality \eqref{e:det} holds for a real number $s\gg
0$ since $Z^e_H(s+n)$ is defined only for $\Re(s)\gg 0$. But the
right side of \eqref{e:det} gives the meromorphic extension over
$\bC$ by Proposition \ref{p:mero-even}, and this also gives the
meromorphic extension of the left side of \eqref{e:det} over
$\bC$.
\end{proof}

\begin{rem}\label{rem:even} \em
Let us remark that $\Det(\Dd^2,s)\neq \Det(\Dd^2,-s)$ as a
meromorphic function over $\bC$. This is because the equality
\eqref{e:det} holds for a real number $s\gg 0$ a priori and the
meromorphic extension of $\Det(\Dd^2,s)$ is given by the right side
of \eqref{e:det}, which does not satisfy this property.
\end{rem}

\vspace{0.3cm}

Now let us show that $Z^e_H(s)$ has the meromorphic extension over
$\bC$. We consider an even smooth function $g(|u|)$ where $g(u)$
is the given function in Section \ref{s6}. We set
$$
H_s(\sigma_{\pm},\lambda)=\int^{\infty}_{-\infty}
g(|u|)e^{-s|u|}e^{i\lambda u}du
$$
for a complex parameter $s$. Then integration by parts gives
\begin{equation}\label{e:eg}
H_s(\sigma_{\pm},\lambda)=\frac{1}{s-i\lambda}\int^{\infty}_0
g'(u)e^{-(s-i\lambda)u}du+\frac{1}{s+i\lambda}\int^{\infty}_{0}
g'(u)e^{-(s+i\lambda)u}du.
\end{equation}
Let us apply the Selberg trace formula to the one parameter family
of functions $f_{s}$ on $G$ with $\Re(s)\gg 0$, such that
$\hat{f}_s(\sigma_{\pm},i\lambda)= H_s(\sigma_{\pm},\lambda)$,
$\hat{f}_s(\sigma,i\lambda)=0$ if $\sigma\neq \sigma_{\pm}$. Then
we get
\begin{align}\label{e:even-sca}
&\sum_{\lambda_j\in\sigma^+_p}H_s(\sigma_+,\lambda_j)
+\sum_{\lambda_j\in\sigma^-_p}H_s(\sigma_-,\lambda_j)\\
&-\frac{\dsi}{4\pi}\int^{\infty}_{-\infty} H_s(\sigma_{+},
\lambda) \ \tr(C_{\Gamma}(\sigma_{+},-i\lambda)
C'_{\Gamma}(\sigma_{+},i\lambda))\
d\lambda \notag \\
=&\sum_{\gamma:\text{hyperbolic}}
l(C_{\gamma})j(\gamma)^{-1}D(\gamma)^{-1}\big(\overline{\chi_{\sigma_{+}}(m_{\gamma})}+
\overline{\chi_{\sigma_{-}}(m_{\gamma})}\big)\, e^{-sl(C_{\gamma})} \notag \\
&-\frac{\kappa d(\sigma_{+})}{2\pi}\int^{\infty}_{-\infty}
H_s(\sigma_{+},\lambda) (\psi(i\lambda+ \frac 12) +\psi(-i\lambda
+\frac 12) )\ d\lambda \notag \\
&+\frac{2}{2\pi}\int^{\infty}_{-\infty} H_s(\sigma_{+}, \lambda)
(P_I(\lambda)+P_U(\lambda))\ d\lambda .\notag
\end{align}
Let us recall the equality
$$
\frac{d}{ds}\log Z^e_H(s+n)=\sum_{\gamma:\text{hyperbolic}}
l(C_{\gamma})j(\gamma)^{-1}D(\gamma)^{-1}\big(\overline{\chi_{\sigma_{+}}(m_{\gamma})}+
\overline{\chi_{\sigma_{-}}(m_{\gamma})}\big)\,
e^{-sl(C_{\gamma})},
$$
and now we investigate the other terms in (\ref{e:even-sca}) to
get the meromorphic extension of $\frac{d}{ds}\log Z^e_H(s+n)$
over $\bC$ and to determine its poles.

\textbf{Discrete eigenvalue term:} \ Using (\ref{e:eg}) as before,
we can see that $H_s(\sigma_{\pm},\lambda_j)$ has a meromorphic
extension over $\bC$ and has the simple poles at $i\lambda_j$ and
$-i\lambda_j$ for $\lambda_j\in\sigma^{\pm}_p$ with the residue
$m_j$ where $m_j$ is the multiplicity of $\lambda_j$.

\textbf{Scattering term:} \ Now we consider the scattering term
$$
-\frac{\dsi}{4\pi}\int^{\infty}_{-\infty}H_s(\sigma_{+},\lambda) \
\tr(C_{\Gamma}(\sigma_{+},-i\lambda)
C'_{\Gamma}(\sigma_{+},i\lambda))\ d\lambda.
$$
As in the previous case, we can show that the function
$$
\Psi(z):=\tr(C_{\Gamma}(\sigma_{+},-z)
\frac{d}{dz}C_{\Gamma}(\sigma_{+},z))
$$
has the following form over $\bC$:
$$
\Psi(z) \ = \ \sum_{\Re(q_k)<0} -\frac{2\,
\Re(q_k)}{(z-q_k)(z+\bar{q}_k)} + \log \ p
$$
for some constant $p$. Here the sum is taken over the set of poles
of $\det C_{\Gamma}(\sigma_+,z)$. We repeat the method in Section
\ref{s6} to prove that the scattering term has a meromorphic
extension over $\bC$ and has poles at $q_k$ for $\Re(q_k)<0$ with
residues $\dsi b_k$ . As before $b_k$ denotes the order of the
pole of $\det C_{\G}(\sigma_+,z)$ at $q_k$.

\textbf{Identity and Unipotent term:}\  It follows from
proposition 3.9 in \cite{GW} that
\begin{equation*}
\frac{1}{2\pi}\int^{\infty}_{-\infty} H_s(\sigma_+,\lambda)
(P_I(\lambda)+P_U(\lambda))\ d\lambda = 0.
\end{equation*}
We now turn our attention to the other unipotent terms. As in
proposition 3.7 of \cite{GW}, we use the Cauchy integral formula
to get
\begin{align*}
\frac{\kappa d(\sigma_{+})}{2\pi}\int^{\infty}_{-\infty}
H_s(\sigma_{+},\lambda) (\psi(i\lambda+ \frac 12) +\psi(-i\lambda
+\frac 12))\ d\lambda \ =\ 2\kappa d(\sigma_+) \psi(s+\frac 12)
\end{align*}
for $\Re(s) \gg 0$ and the right side gives us the meromorphic
extension over $\bC$ of the left side.

Considering the equality (\ref{e:even-sca}) and the above analysis
of other terms, we conclude

\begin{prop}\label{p:mero-even} The Selberg zeta function of
even type $Z^e_H(s)$ has a meromorphic extension over $\bC$.
\end{prop}

\begin{rem}\em
The zeros of $Z^e_H(s+n)\Gamma(s+\frac 12)^{-2\kappa d(\sigma_+)}$
are located at $i\lambda_j$, $-i\lambda_j$ for $\lambda_j\in
\sigma^{\pm}_p$, and at $q_k$ for $\Re(q_k)<0$ and their orders
are $m_j$, $d(\sigma_+)b_k$.
\end{rem}

We now prove the functional equation for $\Det(\Dd^2,s)$ and
$Z^e_H(s)$. First let us define
$$
\mathcal{Z}(s):=\frac{d}{ds}\log Z^e_H(s+n)-2\kappa d(\sigma_{+})
\psi(s+\frac 12)
$$
and we can see that
$$
R(s):=\mathcal{Z}(s)+\mathcal{Z}(-s)+\dsi\Psi(s)
$$
is an even entire function of $s$. Now let $h(s)$ be an even
function which decreases sufficiently rapidly as
$\text{Im}(s)\to\infty$ in the strip $ \{ \ s\in\mathbb{C} \ | \
|\Re(s)| < n+\epsilon, \ \epsilon>0\ \}$. We follow Section
\ref{s6} and consider the contour integral
$$
\mathcal{L}_T:=\frac{1}{2\pi i}\int_{L_T} h(s)\mathcal{Z}(s)\, ds
$$
where $L_T$ is the rectangle with the corners $a+ iT, a-iT, -a+iT,
-a-iT$ with $n<a <n+\epsilon$. Then we have
$$
\lim_{T\to\infty}\mathcal{L}_T =\frac{1}{2\pi
i}\int^{a+i\infty}_{a-i\infty}h(s)(\mathcal{Z}(s)
+\mathcal{Z}(-s))\, ds +\frac{2}{2\pi
i}\int^{-a-i\infty}_{-a+i\infty}h(s) \mathcal{Z}(s)\, ds.
$$
The Cauchy integral theorem gives
\begin{equation}\label{e:eLT}
\lim_{T\to\infty}\mathcal{L}_T= 2\sum_{\lambda_j}m_jh(i\lambda_j)
+\sum_{-a< \Re(q_k)< 0} \dsi b_k h(q_k).
\end{equation}
The simple poles at $q_k$ of $\mathcal{Z}(s)$ in the strip between
$\Re(s)=-a$ and $\Re(s)=0$ have residues $\dsi b_k$ and give the
equality
\begin{equation}\label{e:Eright}
\frac{1}{2\pi i}\int^{-a-i\infty}_{-a+i\infty}h(s)
\mathcal{Z}(s)\, ds = \frac{1}{2\pi
i}\int^{-i\infty}_{i\infty}h(\lambda) \mathcal{Z}(\lambda)\,
d\lambda +\sum_{-a< \Re(q_k)< 0} \dsi b_k h(q_k).
\end{equation}
Similarly, the simple poles at $-q_k$ of $\mathcal{Z}(-s)$ in a
strip between $\Re(s)=a$ and $\Re(s)=0$ have residues $-\dsi b_k$
and provide us with
\begin{align}\label{e:eR(s)}
&\frac{1}{2\pi i}\int^{a+i\infty}_{a-i\infty}h(s)(\mathcal{Z}(s)
 +\mathcal{Z}(-s)) \, ds\\
=&\frac{1}{2\pi
i}\int^{i\infty}_{-i\infty}h(\lambda)(\mathcal{Z}(\lambda)
+\mathcal{Z}(-\lambda))\, d\lambda +\sum_{0< \Re(-q_k)< a}-\dsi
b_k h(-q_k).\notag
\end{align}
We use (\ref{e:eLT}), (\ref{e:Eright}) and (\ref{e:eR(s)}) to
obtain
\begin{align}\label{e:esecond}
&2\sum_{\lambda_j}m_jh(i\lambda_j) = 2\frac{1}{2\pi
i}\int^{-i\infty}_{i\infty}h(\lambda)
(\mathcal{Z}(\lambda)+\frac{\dsi}{2}\Psi(\lambda))\ d\lambda\\
&\qquad\qquad\qquad\qquad + \frac{1}{2\pi
i}\int^{i\infty}_{-i\infty}h(\lambda)(\mathcal{Z}(\lambda)
+\mathcal{Z}(-\lambda)+ \dsi \Psi(\lambda))\ d\lambda.\notag
\end{align}
Replacing $\lambda$ by $i\lambda$ in the first term of the right
side of the equality (\ref{e:esecond}), we see that this term is
equal to
\begin{align}\label{e:eselberg}
&\quad 2 \sum_{\gamma:\text{hyperbolic}}
l(C_{\gamma})j(\gamma)^{-1}D(\gamma)^{-1}\big(\overline{\chi_{\sigma_{+}}(m_{\gamma})}+
\overline{\chi_{\sigma_{-}}(m_{\gamma})}\big)\ \frac{1}{2\pi}
\int^{\infty}_{-\infty} h(i\lambda)
e^{-i\lambda l(C_{\gamma})} \ d\lambda \\
&-2\frac{1}{2\pi}\int^{\infty}_{-\infty}h(i\lambda)(2\kappa
d(\sigma_{+}) \psi(i\lambda+\frac 12) ) \ d\lambda
+2\frac{\dsi}{4\pi }\int^{\infty}_{-\infty}h(i\lambda)
\Psi(i\lambda) \ d\lambda.\notag
\end{align}
By (\ref{e:esecond}), (\ref{e:eselberg}) and the Selberg trace
formula applied to $h(i\lambda)$, we have
$$
\frac{1}{2\pi i}\int^{i\infty}_{-i\infty}h(\lambda)(R(\lambda)+
4P_I(-i\lambda)+4P_U(-i\lambda))\ d\lambda =0
$$
for any even holomorphic $h(s)$ satisfying suitable growth
conditions. Since
$$
R(s)+4P_I(-is)+4P_U(-is)=R(s)+4P_I(is)+4P_U(is)
$$
is an entire function, we have
$$
\mathcal{Z}(s) +\mathcal{Z}(-s)+4P_I(is)+4P_U(is)=- \dsi \Psi(s)
$$
over $\bC$. Now we get the equality
$$
\frac{d}{ds}\log \Bigr( {Z^e_H}(n+s)\Gamma(s+\frac 12)^{-2\kappa
d(\sigma_{+})} \Bigr)+4P_I(is)+4P_U(is)
$$
$$=\frac{d}{ds}\log\Bigr({Z^e_H}(n-s)\Gamma(-s+\frac 12)^{-2\kappa
d(\sigma_{+})} \Bigr) - \dsi\frac{d}{ds} \log(\det
C_{\Gamma}(\sigma_{+},s)),
$$
which leads to the formula
\begin{align}\label{e:Efe}
&Z^e_H(n+s)\Gamma(s+\frac 12)^{-2\kappa d(\sigma_{+})}
\exp(4\int^s_0
P_I(i\lambda)+P_U(i\lambda)\ d\lambda)\\
=& \  ( \det C_{\Gamma}(\sigma_+,s))^{-\dsi}\
Z^e_H(n-s)\Gamma(-s+\frac 12)^{-2\kappa d(\sigma_{+})}  .\notag
\end{align}
This equality holds a priori up to constant $c$. As in the proof
of lemma 4.3 in \cite{GW}, we multiply $s^{-2m_0}$ to both sides
of (\ref{e:Efe}) where $m_0$ is the multiplicity of the zero
eigenvalue of $\Dd$. Doing this removes the spectral zero of
$Z^e_H(n\pm s)$ at $s=0$. If we compare the remaining parts at
$s=0$, we can see that the constant $c$ equals $1$. From
(\ref{e:det}), we have
\begin{align*}
\Det(\Dd^2,s)^2 \ = \ & C^2\, Z^e_H(s+n)^2 \Gamma(s+\frac
12)^{-4\kappa d(\sigma_{+})} \exp(4\int^s_0
P_I(i\lambda)+P_U(i\lambda)\ d\lambda).
\end{align*}
Finally, we combine (\ref{e:Efe}) and this equality to get

\begin{thm}\label{t:fd} For any $s\in\mathbb{C}$, we have
\begin{multline}\label{e:final}
\Det(\Dd^2,s)^2 \\ = C^2\, ( \det C_{\Gamma}(\sigma_+,s))^{-\dsi}\
Z^e_H(n+s)Z^e_H(n-s) \Bigr(\Gamma(s+\frac 12)\Gamma(-s+\frac
12)\Bigr)^{-2\kappa d(\sigma_{+})}.
\end{multline}
\end{thm}

\begin{rem}\em The right side of (\ref{e:final}) has the zeros at $i\lambda_j$, $-i\lambda_j$ of order $2m_j$
for $\lambda_j\in\sigma^{\pm}_p$,
 and at $q_k$ of order $2d(\sigma_+)b_k$ for $\Re(q_k)<0$ . These are the zeros
of $\Det(\Dd^2,s)^2$ as we expected.
\end{rem}

\end{document}